\documentclass[article, onesided, a4paper]{amsart}

\usepackage[utf8]{inputenc}

\usepackage{amsmath}
\usepackage{amssymb}

\usepackage{bm}

\usepackage{physics}

\usepackage{xurl}

\usepackage{imakeidx}

\usepackage{hyperref}

\let\oldproof\proof
\def\proof{\oldproof\unskip}

\usepackage{etoolbox}

\makeatletter

\patchcmd\deferred@thm@head
  {\addvspace{-\parskip}}
  {}
  {}{\typeout{\string\deferred@thm@head patch failed!}}

\makeatletter 

\usepackage{mathtools}
\mathtoolsset{showonlyrefs}
\usepackage{enumitem}
\usepackage{amsthm}

\usepackage{tikz-cd}

\newtheorem{thm}[subsection]{Theorem}
\newtheorem{lem}[subsection]{Lemma}
\newtheorem{prop}[subsection]{Proposition}
\newtheorem{cor}[subsection]{Corollary}
\newtheorem{fact}[subsection]{Fact}

\theoremstyle{remark}
\newtheorem*{rmk}{Remark}

\newtheorem*{note}{Note}

\theoremstyle{definition}
\newtheorem{defn}[subsection]{Definition}
\newtheorem*{nota}{Notation}
\newtheorem{eg}[subsection]{Example}
\newtheorem{egs}[subsection]{Examples}

\newcommand{\tp}{\mathrm{tp}}

\DeclareMathOperator{\ch}{char}

\DeclareMathOperator{\Th}{Th}
\newcommand{\id}{\mathrm{id}}

\DeclareMathOperator{\ord}{ord}

\DeclareMathOperator{\trdeg}{trdeg}

\newcommand{\alg}{\mathrm{alg}}

\newcommand{\vf}{\mathrm{vf}}
\newcommand{\vd}{\mathrm{vd}}
\newcommand{\qf}{\mathrm{qf}}

\DeclareMathOperator{\Jet}{Jet}

\DeclareMathOperator{\supp}{supp}
\DeclareMathOperator{\res}{res}

\newcommand{\Reg}{\mathrm{Reg}}

\newcommand{\DL}{\mathrm{DL}}

\newcommand{\bdelta}{\bm{\delta}}

\newcommand{\Hom}{\mathrm{Hom}}

\newcommand{\ac}{\mathrm{ac}}

\renewcommand{\leq}{\leqslant}
\renewcommand{\geq}{\geqslant}
\renewcommand{\preceq}{\preccurlyeq}

\renewcommand{\epsilon}{\varepsilon}
\renewcommand{\O}{\mathcal{O}}
\renewcommand{\phi}{\varphi}
\newcommand{\m}{\mathfrak{m}}

\newcommand{\la}{\langle}
\newcommand{\ra}{\rangle}
\newcommand{\eqra}{{\equiv\ra}}

\renewcommand{\L}{\mathcal{L}}
\newcommand{\Lvf}{\mathcal{L}_\mathrm{vf}}
\newcommand{\Lvd}{\mathcal{L}_\mathrm{vd}}

\newcommand{\Mor}{\mathrm{Mor}}

\newcommand{\Pas}{\mathrm{Pas}}

\newcommand{\UC}{\mathrm{UC}}

\DeclareMathOperator{\DCF}{DCF}

\newcommand{\DCVF}{\mathrm{DCVF}}
\newcommand{\ACVF}{\mathrm{ACVF}}

\newcommand{\DHen}{\mathrm{DH}}

\newcommand{\Z}{\mathbb{Z}}

\newcommand{\Q}{\mathbb{Q}}
\newcommand{\A}{\mathbb{A}}

\newcommand{\rest}{{\upharpoonright}}

\newcommand{\ddt}{\frac{\dd}{\dd t}}

\newcommand{\AKE}{\mathrm{AKE}}

\newcommand{\Alg}{\mathbf{\text-Alg}}

\def\Ind#1#2{#1\setbox0=\hbox{$#1x$}\kern\wd0\hbox to 0pt{\hss$#1\mid$\hss}
\lower.9\ht0\hbox to 0pt{\hss$#1\smile$\hss}\kern\wd0}

\def\notind#1#2{#1\setbox0=\hbox{$#1x$}\kern\wd0
\hbox to 0pt{\mathchardef\nn=12854\hss$#1\nn$\kern1.4\wd0\hss}
\hbox to 0pt{\hss$#1\mid$\hss}\lower.9\ht0 \hbox to 0pt{\hss$#1\smile$\hss}\kern\wd0}

\usepackage{graphicx}

\makeatletter
\renewenvironment{proof}[1][\proofname]{\par
  \vspace{-\topsep}% remove the space after the theorem
  \pushQED{\qed}%
  \normalfont
  \topsep0pt \partopsep0pt % no space before
  \trivlist
  \item[\hskip\labelsep
        \itshape
    #1\@addpunct{.}]\ignorespaces
}{%
  \popQED\endtrivlist\@endpefalse
  \addvspace{6pt plus 6pt} % some space after
}
\makeatother

\title{Differentially Henselian Fields}
\author{Gabriel Ng}
\date{10 February 2025}

\begin{document}

\begin{abstract}
We study the class of differentially henselian fields, which are henselian valued fields equipped with generic derivations in the sense of Cubides Kovacics and Point, and are special cases of differentially large fields in the sense of Le\'on S\'anchez and Tressl. We prove that many results from henselian valued fields as well as differentially large fields can be lifted to the differentially henselian setting, for instance Ax-Kochen/Ershov principles, characterisations in terms of differential algebras, etc. We also give methods to concretely construct such fields in terms of iterated power series expansions and inductive constructions on transcendence bases.
\end{abstract}

\maketitle

\section{Introduction}

In \cite{SanchezTressl2020} and \cite{SanchezTressl2023}, Le\'on-S\'anchez and Tressl introduce and study the class of \emph{differentially large fields}, which are large fields equipped with a derivation which satisfies a certain existential closure condition, i.e. a differential field $(K, \delta)$ is differentially large if and only if for every differential field extension $(L, \partial)$ such that $K$ is existentially closed in $L$, $(K, \delta)$ is existentially closed in $(L, \partial)$. In this paper, we consider a similar notion for henselian valued fields: we say that a valued-differential field $(K, v, \delta)$ is \emph{differentially henselian} if it is (nontrivially) henselian as a pure valued field, and for any valued-differential field extension $(L, w, \partial)$ of $(K, v, \delta)$, if $(K, v)$ is existentially closed in $(L, w)$ as pure valued fields, then $(K, v, \delta)$ is existentially closed in $(L, w, \partial)$ as valued-differential fields (Later Definition \ref{dh_field_defn}). We will only consider the case of characteristic 0 fields equipped with a single derivation.

Differentially henselian fields are examples of differential fields with a \emph{generic derivation} as studied by Guzy, Cubides Kovacsics and Point in \cite{GuzyPoint2010} and \cite{CubidesPoint2019}, i.e. they are precisely the henselian valued fields equipped with a derivation satisfying the axiom scheme (DL) (Theorem \ref{dh_equiv_model_DH_thm}). We also show that they coincide with the models of the theory ($\UC'$) as introduced by Guzy in \cite{Guzy2006} as an valued-differential field analogue for the uniform companion theory (UC) constructed by Tressl for pure differential fields.

In this paper, we will establish a wide variety of results concerning the class of differentially henselian fields. In Section \ref{section_existential_lifting}, we show that a certain `existential lifting' property can be lifted from the level of valued fields to valued-differential fields. An important byproduct of this is the Relative Embedding Theorem (\ref{relative_embedding_thm}), which we use to establish a number of results in the following sections.

In particular, in Section \ref{section_dh_ake}, we show how certain Ax-Kochen/Ershov principles can be lifted from the setting of henselian valued fields to differentially henselian fields. That is, an elementary class of differentially henselian fields satisfies the AKE principles for existential closure, completeness, model completeness, and relative subcompleteness if and only if the corresponding class of pure valued fields also satisfies the same property (Propositions \ref{ake_exists_prop}, \ref{ake_equiv_prop} and Theorem \ref{ake_prec_thm}). In Section \ref{section_dh_qe_for_ac}, we prove that differentially henselian fields of equicharacteristic 0 equipped with an angular component map eliminates field sort quantifiers in the differential Pas language (Theorem \ref{dh_ac_has_qe}). In Section \ref{section_dh_stable_embedded}, we show that in equicharacteristic 0 and unramified mixed characteristic differentially henselian fields, the residue field and value group are stably embedded (Theorem \ref{dh_res_field_val_gp_stabl_emb_thm}), and the induced structure on the constants is precisely that induced by the dense pair of valued fields (Theorem \ref{induced_str_constants_is_pair_str_thm}).

In Section \ref{section_constructions}, we demonstrate methods through which differentially henselian fields may be constructed. Of particular interest is a construction by iterating Hahn series extensions (Theorem \ref{K_infty_dh}), which gives a natural differentially henselian extension of any henselian valued-differential field. We also show by an inductive argument that any henselian valued field of infinite transcendence degree admits a derivation such that the resulting valued-differential field is differentially henselian (Corollary \ref{inf_trdeg_diff_hensel}). In Section \ref{section_sections}, we construct differential sections of the residue field of a differentially henselian field equipped with an arbitrary derivation.

In the latter sections, we study the connections between differentially large and differentially large henselian fields. We apply a theorem of Widawski to show that every non-algebraically closed differentially large field with a henselian valuation is differentially henselian (Proposition \ref{diff_large_not_ac_implies_diff_hensel_prop}). In Section \ref{section_equiv_chars}, we prove a variety of equivalent characterisations of differentially henselian fields in the style of \cite[Theorem 4.3]{SanchezTressl2020}. In particular, we show that a valued-differential field $(K, v, \delta)$ is differentially henselian if and only if it is existentially closed in the field of formal Laurent series $(K((t)), v\circ v_t, \delta + \ddt)$ (Theorem \ref{dh_iff_ec_laurent}). We also show that $(K, v, \delta)$ is differentially henselian if and only if it is henselian as a valued field, and any differentially finitely generated $K$-algebra with a $K$-point has a differential $K$-point `arbitrarily close' with respect to the valuation topology (Theorem \ref{dh_iff_4.3iv(h)}).

The final section of this paper is dedicated to showing that the differential Weil descent as established in \cite{SanchezTressl2018} interacts well with valuations (Theorem \ref{weil_descent_continuous}), and we show that in a restricted setting the converse also holds (Proposition \ref{weil_descent_continuous_converse}). We apply this new tool to show that any algebraic extension of a differentially henselian field (equipped with the unique extensions of the derivation and valuation) is again differentially henselian (Theorem \ref{alg_extn_of_dh_is_dh}).

\section{Preliminaries}

Throughout this article, every field is of characteristic 0, and every ring will be commutative and unital.

Let $K$ be a ring. Recall that a \textbf{derivation} on $K$ is an additive map $\delta: K \to K$ satisfying the product rule, i.e. for any $a, b \in K$,
\[
\delta(ab) = \delta(a)b + a \delta(b).
\]
A \textbf{differential ring/field} is a ring/field $K$ equipped with a \emph{single} derivation $\delta: K \to K$. 
A \textbf{differential subring} of $(K, \bdelta)$ is a subring $L$ of $K$ such that $L$ is closed under $\delta$. A \textbf{differential field} is a differential ring which is a field. A \textbf{differential ideal} of a differential ring $(A, \delta)$ is an ideal $\mathfrak{p}$ of $A$ such that $\delta(\mathfrak{p}) \subseteq \mathfrak{p}$.

\begin{nota}
    Let $(K, \delta)$ be a differential ring, and $a \in K$. For $n < \omega$, the \textbf{$n$-jet} of $a$ is the tuple $\Jet_n(a) = (a, a',..., a^{(n)})$.
    The (infinite) \textbf{jet} of $a$ is the $\omega$-tuple $\Jet(a) = (a, a', a'',...)$. For a tuple $\bar{a} \in K$, write $\Jet_n(\bar{a})$ for the concatenation of the $n$-jets of each entry of $\bar{a}$, and similarly $\Jet(\bar{a})$.
\end{nota}

A \textbf{differential polynomial} in variables $\bar{x} = (x_i)_{i < \alpha}$ over a differential ring $(K, \delta)$ is a polynomial with variables among $x_i^{(n)}$, where $i<\alpha$, $n<\omega$. The \textbf{ring of differential polynomials} (over $K$ in the variables $\bar{x}$) is the polynomial ring $K[x_i^{(n)} : i<\alpha, n<\omega]$ equipped with the standard derivation $\delta$ extending the derivation on $K$ satisfying $\delta(x_i^{(n)}) = x_i^{(n+1)}$.

For a differential polynomial (in one variable) $f \in K\{x\}$, the \textbf{order of $f$}, denoted $\ord(f)$ is the largest $n$ such that $x^{(n)}$ occurs in $f$. If $\ord(f) = n$, \textbf{separant of $f$} is the formal partial derivative $s(f) = \frac{\partial f}{\partial x^{(n)}}$.

\begin{nota}
    For a differential polynomial $f(\bar{x}) \in K\{\bar{x}\}$, we write $f_\alg(\bar{x})$ for the (non-differential) polynomial obtained by considering $f$ as a non-differential polynomial by forgetting the derivation.
\end{nota}

Let $(K, \delta)$ be a differential ring. A \textbf{differential $K$-algebra} is a $K$-algebra $A$ equipped with a derivation $\partial$ such that the structure map $\eta_A: K \to A$ is a differential ring homomorphism $(K, \delta) \to (A, \partial)$. A differential $K$-algebra $A$ is said to be \textbf{differentially finitely generated} if there is a surjective differential $K$-algebra homomorphism $K\{\bar{x}\} \to A$, where $\bar{x}$ is a tuple of finitely many indeterminates. If $K$ is a differential field, and $\bar{a}$ is a tuple of elements from some differential field extension of $K$, $K\la \bar{a} \ra$ denotes the \textbf{differential field generated by $\bar{a}$ over $K$}.

For $K$ a (differential) field and $A$ a (differential) $K$-algebra, a (differential) \textbf{$K$-point of $A$} is a (differential) $K$-algebra homomorphism $\phi: A \to K$.

Recall that a field $K$ is said to be \textbf{large} (also known as \textbf{ample}) if $K$ is existentially closed in the field of Laurent series $K((t))$. Examples include algebraically closed, real closed, and fields admitting a nontrivial henselian valuation.

\begin{defn}[{\cite[Definition 4.2]{SanchezTressl2020}}]
    A differential field $(K, \delta)$ is \textbf{differentially large} if $K$ is large as a field, and for any differential field extension $(L, \partial)$ of $(K, \delta)$, if $K$ is existentially closed in $L$ as a pure field, then $(K, \delta)$ is existentially closed in $(L, \partial)$ as a differential field.
\end{defn}

We recall a number of equivalent characterisations of differential largeness:
\begin{thm} \label{diff_large_eq_conds_thm}
    Let $(K, \delta)$ be a differential field. The following are equivalent:
    \begin{enumerate}
        \item $(K, \delta)$ is differentially large;
        \item $(K, \delta)$ is existentially closed in $(K((t)), \delta + \ddt)$;
        \item $K$ is large as a field, and any differentially finitely generated $K$-algebra with a $K$-point $\phi: A \to K$ has a differential $K$-point $\psi: A \to K$;
        \item $K$ is large as a field, and for any pair of differential polynomials $f, g \in K\{x\}$ with $\ord(f) \geq \ord(g)$, if there is $\bar{a} \in K$ such that $f_\alg(\bar{a}) = 0$, $s(f)_\alg(\bar{a}) \neq 0$ and $g(\bar{a}) \neq 0$, then there is $b \in K$ with $f(b) = 0$ and $g(b) \neq 0$.
    \end{enumerate}
\end{thm}

A \textbf{valuation} on a field $K$ is a surjective map $v: K \to \Gamma \cup \{\infty\}$, where $\Gamma$ is an ordered abelian group, satisfying, for any $a, b \in K$:
\begin{itemize}
    \item $v(a) = \infty$ if and only if $a = 0$;
    \item $v(ab) = v(a) + v(b)$ (where $\gamma + \infty = \infty = \infty + \gamma$ for any $\gamma \in \Gamma \cup \{\infty\}$);
    \item $v(a + b) \geq \min(v(a), v(b))$ (where $\infty > \gamma$ for all $\gamma \in \Gamma$).
\end{itemize}
A \textbf{valued field} is a field equipped with a valuation $v$, denoted $(K, v)$. The \textbf{value group} of $(K, v)$ is the ordered abelian group $v(K^\times)$, denoted $vK$. The \textbf{valuation ring} of $(K, v)$ is the subring $\O_v = \{a \in K: v(a) \geq 0\}$ consisting of elements of nonnegative valuation. The valuation ring is a local ring with maximal ideal $\m_v = \{a \in K : v(a) > 0\}$, and the quotient $\O_v/\m_v$ is called the \textbf{residue field}, denoted $Kv$. The quotient map $\O_v \to Kv$ is denoted $\res$ (or $\res_v$ if the valuation needs to be specified). The \textbf{characteristic} of a valued field $(K, v)$ is the pair $(\ch K, \ch Kv)$. If $(K, v)$ has characteristic $(0, 0)$ or $(0, p)$ where $p$ is positive, we say that $(K, v)$ is \textbf{equicharacteristic 0} or \textbf{mixed characteristic}, respectively.

For a valued field $(K, v)$ and a valuation $w$ on the residue field $Kv$, we denote the \textbf{composition of $w$ with $v$} by $w \circ v$, i.e. this is the valuation on $K$ corresponding to the valuation ring $\res_v^{-1}(\O_w)$.

A valued field $(K, v)$ is said to be \textbf{henselian} it satisfies \textbf{Hensel's Lemma}:
\begin{center}
    Let $f \in \O_v[x]$, and $a \in \O_v$ such that $\res(f(a)) = 0$ and $\res(f'(a)) \neq 0$, where $f'$ denotes the formal derivative of $f$. Then, there is $b \in \O_v$ such that $f(b) = 0$, and $\res(a - b) = 0$.
\end{center}
Equivalently, $(K, v)$ is henselian if $v$ admits a unique extension to any algebraic extension of $K$.
We assume throughout this article that henselian valuations are \textbf{nontrivial}, i.e. the valuation has nontrivial value group. A pure field $K$ is said to be \textbf{henselian} if there exists a valuation $v$ on $K$ such that $(K, v)$ is henselian.

\begin{fact}
    Every field which admits a henselian valuation is large.
\end{fact}

A crucial fact which we will apply liberally is that henselian valued fields satisfy a version of the \emph{implicit function theorem} \cite{Kuhlmann2000}:

\begin{thm}[Implicit Function Theorem for Henselian Fields] \label{henselian_ift}
    Let $(K, v)$ be a henselian valued field. Let $f_1,...,f_n \in K[x_1,...,x_m,y_1,...,y_n]$. Set 
        \[
        \bar{z} = (x_1,...,x_m,y_1,...,y_n)
        \]
        and
        \[
        J(\bar{z}) = \begin{pmatrix}
        \frac{\partial f_1}{\partial y_1}(\bar{z}) & \cdots & \frac{\partial f_1}{\partial y_n} (\bar{z}) \\
        \vdots & \ddots & \vdots \\
        \frac{\partial f_n}{\partial y_1}(\bar{z}) & \cdots & \frac{\partial f_n}{\partial y_n}(\bar{z})
        \end{pmatrix}.
        \]
        Suppose that $\bar{f} = (f_1,...f_n)$ has a zero $(\bar{a}, \bar{b}) = (a_1,...,a_m,b_1,...,b_n) \in K^{n+m}$, and the determinant of $J(\bar{a}, \bar{b})$ is nonzero. Then, there is some $\alpha \in vK$ such that for all $\bar{c} = (c_1,...,c_m) \in K^m$ with $v(a_i - c_i) > 2\alpha$ for each $i$, there is a unique $\bar{d} = (d_1,...,d_n) \in K^n$ such that $(\bar{c}, \bar{d})$ is a zero of $\bar{f}$, and $v(b_i - d_i) > \alpha$ for each $i$.
\end{thm}

\section{Differential Henselianity}

\begin{defn} \label{dh_field_defn}
    A \textbf{valued-differential field} is a field $K$ equipped with a valuation $v$ and derivation $\delta$, denoted by the ordered triple $(K, v, \delta)$. No interaction is specified between $v$ and $\delta$.
\end{defn}

\begin{rmk}
    In particular, we do not assume that the derivation is \emph{small} or continuous with respect to the valuation topology. In fact, most derivations we consider will be highly discontinuous.
\end{rmk}

\begin{defn}
    The \textbf{language of valued fields}, denoted $\Lvf$, consists of the language of rings along with a binary relation symbol $|$ interpreted as `valuation divisibility', i.e. for a valued field $(K, v)$ and $a, b \in K$, $a|b$ holds if and only if $v(a) \leq v(b)$. The \textbf{language of valued-differential fields}, denoted $\Lvd$, is the language of valued fields along with a unary function symbol $\delta$ interpreted as the derivation.
\end{defn}

We define the notion of \emph{differential henselianity} analogously to differential largeness in the context of valued-differential fields.

\begin{defn}
    We say that a valued-differential field $(K, v, \delta)$ is \textbf{differentially henselian} if $(K, v)$ is henselian as a pure valued field, and for any valued-differential field extension $(L, w, \partial)$ of $(K, v, \delta)$, if $(K, v)$ is existentially closed in $(L, w)$ as a pure valued field, then $(K, v, \delta)$ is existentially closed in $(L, w, \partial)$ as valued-differential fields.
\end{defn}

\begin{lem} \label{diff_hens_implies_diff_large_lem}
    Every differentially henselian field is differentially large as a differential field.
\end{lem}
\begin{proof}
    Let $(K, v, \delta)$ be differentially henselian. Let $(L, \partial)$ be a differential field extension of $(K, \delta)$ such that $K$ is existentially closed in $L$ as a pure field. By existential closure, $L$ embeds as a field over $K$ in any sufficiently saturated elementary extension of $K$. Let $(K^*, v^*, \delta^*)$ be such an extension of $(K, v, \delta)$ and $\phi: L \to K^*$ be an embedding of fields over $K$. Let $w$ be the valuation on $L$ obtained by restricting $v^*$ to the image of $L$ under $\phi$. By construction, $(L, w)$ is a valued field extension of $(K, v)$ which embeds in an elementary extension of $(K, v)$, thus $(K, v)$ is existentially closed in $(L, w)$. By differential henselianity of $(K, v, \delta)$, we obtain that $(K, v, \delta)$ is existentially closed in $(L, w, \partial)$ as valued-differential fields. Taking reducts, we have that $(K, \delta)$ is existentially closed in $(L, \partial)$, as required.
\end{proof}

\section{Differential Lifting} \label{section_existential_lifting}

In this section, we will axiomatise the class of differentially henselian fields in terms of a differential lifting condition, as well as prove a useful relative embedding theorem which we shall exploit in later sections. We will prove the following characterisation of differential henselianity:

\begin{thm} \label{dh_equiv_model_DH_thm}
    Let $(K, v, \delta)$ be a valued-differential field, henselian as a pure valued field. Then, $(K, v, \delta)$ is differentially henselian if and only if $(K, v, \delta)$ satisfies the following axiom scheme:
    \begin{center}
            For any differential polynomial $f(x) \in K\{x\}$ of order $n$ and any $\bar{a} \in K^{n+1}$ such that $f_\alg(\bar{a}) = 0$ and $s(f)_\alg(\bar{a}) \neq 0$, there is, for any $\gamma \in vK$, $b \in K$ such that $f(b) = 0$ and $\Jet_n(b) \in B_\gamma(\bar{a})$. 
    \end{center}
\end{thm}

\begin{rmk}
    This says that the class of differentially henselian fields are precisely those valued-differential fields which are henselian and satisfy the axiom scheme (DL) from \cite{CubidesPoint2019}. Following Cubides Kovacsics and Point, we refer to this axiom scheme in the language $\Lvd$ as (DL).
\end{rmk}

\begin{defn}
    Denote by $\DHen$ the $\Lvd$-theory consisting of the axioms for a henselian valued field of characteristic 0, along with the axiom scheme $(\DL)$ above.
\end{defn}

\begin{lem} \label{dh_implies_model_DH_lem}
    Let $(K, v, \delta)$ be a differentially henselian field. Then, $(K, v, \delta) \models \DHen$.
\end{lem}
\begin{proof}
    By \cite[Theorem 2.3.4]{CubidesPoint2019}, there is an elementary extension $(L, w)$ of $(K, v)$ and a derivation $\partial$ on $L$ extending $\delta$ on $K$ such that $(L, w, \partial) \models \DHen$. By differential henselianity, $(K, v, \delta)$ is existentially closed in $(L, w, \partial)$. Since $(\DL)$ is inductive, by existential closure, we also have $(K, v, \delta) \models (\DL)$, as required.
\end{proof}

For the converse, we will show an `existential lifting' result for models of $\DHen$.

\begin{lem} \label{valuation_atomics}
Let $(K, v)$ be a valued field. Let $\phi(\bar{x})$ be a consistent $\Lvf(K)$ formula which is a boolean combination of valuation atomic formulae. Let $\bar{a} \in K$ be such that $K \models \phi(\bar{a})$. Then, there exists $\gamma \in vK$ such that for any $\bar{b} \in B_\gamma(\bar{a})$, $K \models \phi(\bar{b})$. That is, the set defined by $\phi(\bar{x})$ is open in $K^n$ with respect to the valuation topology. 
\end{lem}
\begin{proof}
% As polynomials are continuous with respect to the valuation topology, for every polynomial $f$ appearing in $\phi$, there is some $\gamma \in vK$ such that for all $\bar{b} \in B_\gamma(\bar{a})$, $v(f(\bar{b}) - f(\bar{a})) > v(f(\bar{a}))$. In particular, for any $\bar{b} \in B_\gamma(\bar{a})$, $v(f(\bar{a})) = v(f(\bar{b}))$. Fix $\gamma \in vK$ such that the above holds for all polynomials $f$ appearing in $\phi$. Now, for any atomic formula $\psi(\bar{x})$ of the form $v(f(\bar{x})) \leq v(g(\bar{x}))$ appearing in $\phi$, and for any $b \in B_\gamma(\bar{a})$, $K \models \psi(\bar{a})$ if and only if $K \models \psi(\bar{b})$. As the truth of every atomic formula is preserved by replacing $\bar{a}$ with $\bar{b}$, we also have that $K \models \phi(\bar{b})$.
First note that the set $\{(x, y) \in K^2: v(x) \leq v(y)\}$ is clopen in $K^2$, and for any polynomials $f(\bar{x}), g(\bar{x})$ over $K$, the map $K^n \to K^2: \bar{x} \mapsto (f(\bar{x}), g(\bar{x}))$ is continuous. Thus, the preimage $\{\bar{x} \in K^n: v(f(\bar{x})) \leq v(g(\bar{x}))\}$ is also clopen. Since the set defined by $\phi(\bar{x})$ is a finite boolean combination of clopen sets, it is also clopen (and hence open).
\end{proof}

\begin{lem}\label{alg_geom_formulae}
Let $(K, v, \delta)$ be a valued-differential field, and $A$ a differential subfield of $K$. Let $\phi(\bar{x})$ be a quantifier-free $\Lvd(A)$-formula, and $\bar{a} \in K$ such that $K \models \phi(\bar{a})$. There are $\Lvd(A)$-formulae $\psi(\bar{x}), \chi(\bar{x})$ such that
\[
\models \forall\bar{x}((\psi(\bar{x}) \wedge \chi(\bar{x})) \to \phi(\bar{x})),
\]
and
\[
K \models \psi(\bar{a})\wedge\chi(\bar{a}),
\]
and $\psi, \chi$ are conjunctions of algebraic and valuation atomic formulae and negations of algebraic and valuation atomic formulae, respectively.
\end{lem}
\begin{proof}
Let $\bigvee_i \phi_i(\bar{x})$ be the disjunctive normal form of $\phi$, and let $\phi_N(\bar{x})$ be the disjunct which holds for $\bar{a}$. Let $\psi(\bar{x})$ and $\chi(\bar{x})$ be the conjunctions of the algebraic and valuation atomic formulae appearing in $\phi_N$, respectively.
\end{proof}

For later convenience, we will call the formulae $\psi$ and $\chi$ in the above lemma the \textbf{algebraic} and \textbf{valuation parts} of $\phi$ (with respect to $\bar{a}$), respectively.

\begin{nota}
    For a valued-differential field $(K, v, \delta)$ and a set of parameters $A \subseteq K$, the quantifier-free $\Lvf$ and $\Lvd$-types of $\bar{a} \in K$ over $A$ are denoted $\qf\tp_\vf(\bar{a}/A)$ and $\qf\tp_\vd(\bar{a}/A)$, respectively.
\end{nota}

We first consider the case where the tuple $\bar{a}$ is differentially algebraically independent over $A$. For the rest of this section, we fix (unless otherwise stated) $(K, v, \delta)$ and $(L, w, \partial)$ valued-differential fields with $(L, w, \partial) \models \DHen$, and $A$ a common valued-differential subfield. We also assume that $(L, w, \partial)$ is $|A|^+$-saturated.

\begin{lem}\label{embedding_transcendentals}
Let $\bar{a} = (a_i)_{i<n} \in K$ be differentially algebraically independent over $A$. Let $\bar{c} = (c_{ij})_{i<n, j<\omega} \in L$ be a realisation of $\qf\tp_\vf(\Jet(\bar{a})/A)$ in $L$. Then, for any $\epsilon \in wL$, there is $\bar{u} \in L$ realising $\qf\tp_\vd(\bar{a}/A)$ and $\Jet(\bar{u}) \in B_\epsilon(\bar{c})$.
\end{lem}
\begin{proof}
We show that the type $\qf\tp_\vd(\bar{a}/A)$ along with `$\Jet(\bar{y}) \in B_\epsilon(\bar{c})$' is finitely satisfiable in $L$. Let $\phi(\bar{y}) \in \qf\tp_\vd(\bar{a}/A)$, and let $\psi$ and $\chi$ be the algebraic and valuation parts of $\phi$, respectively. Let $N$ be the highest order derivative of any variable appearing in $\phi$. Since $\bar{a}$ is differentially algebraically independent over $A$, we may assume that $\psi$ consists only of differential polynomial inequations. As polynomials are continuous with respect to the valuation topology, there is some $\mu \in wL$ such that for any differential polynomial inequation $g(\bar{y}) \neq 0$ appearing in $\psi$, and $nN$-tuple $\bar{u} \in B_\mu((c_{ij})_{i<n, j\leq N})$, we have that $g(\bar{u}) \neq 0$.

Further, by Lemma \ref{valuation_atomics}, there is some $\lambda \in wL$ such that for any $nN$-tuple $\bar{u} \in B_\lambda((c_{ij})_{i<n, j\leq N})$, $\psi(\bar{u})$ holds. Let $\epsilon' = \max(\epsilon, \mu, \lambda)$. By differential henselianity, there is $\bar{u} = (u_i)_{i<n} \in L$ such that $\Jet_N(v_i) \in B_{\epsilon'}((c_{ij})_{j\leq N})$. Thus, we find that $L \models \phi(\bar{v}) \wedge \Jet_N(\bar{v}) \in B_\epsilon((c_{ij})_{i<n, j\leq N})$, as required. Now, apply saturation.
\end{proof}

\begin{lem}\label{approximating_types_restricted_form}
% Let $(K, v, \delta)$ and $(L, w, \partial)$ be valued-differential fields, with $(L, w, \partial) \models \DHen$, and $A$ be a common valued-differential subfield. Suppose that $(L, w, \partial)$ is $|A|^+$-saturated. 
Let $\bar{a}b \in K$ be such that $\bar{a} = (a_i)_{i<n}$ is differentially algebraically independent over $A$, and $b$ is differentially algebraic over $A\la \bar{a}\ra$. Suppose there is $\bar{c}\bar{d} = (c_{ij})_{i<n, j<\omega}(d_k)_{k<\omega} \in L$ realising $\qf\tp_\vf(\Jet(\bar{a})\Jet(b)/A)$. Then, for any $\gamma \in wL$ there is $\bar{\alpha}\beta \in L$ such that $\bar\alpha\beta \models \qf\tp_\vd(\bar{a}b/A)$ and $\Jet(\alpha)\Jet(\beta)\in B_\gamma(\bar{c}\bar{d})$.
\end{lem}
\begin{proof}
Let $f(x)$ be the differential minimal polynomial of $b$ over $A\la \bar{a}\ra$. Let $\ord(f) = m$. By clearing denominators, we may write $f(x)$ as $F(\bar{a}, x)$, where $F(\bar{y}, x)$ is a differential polynomial over $A$ in $n+1$ variables. Denote its separant (with respect to $x$) by $s(F)(\bar{y}, x)$. 

We claim that $p(\bar{y}, x) \coloneqq \qf\tp_\vd(\bar{a}b/A)$ along with the partial type expressing that `$\Jet(\bar{y})\Jet(x) \in B_\gamma(\bar{c}\bar{d})$' is finitely satisfiable in $L$. Let $\phi(\bar{y}, x)$ be some formula in $p$. Let $\psi$ and $\chi$ denote the algebraic and valuation parts of $\phi$, respectively. 

For a tuple $\bar{t}u$ to satisfy $\psi(\bar{y}, x)$, it suffices for $\bar{t}$ to be differentially algebraically independent over $A$, $F(\bar{t}, u) = 0$ and for each differential polynomial inequation $g(\bar{y}, x) \neq 0$ appearing in $\psi$, $g(\bar{t}, u) \neq 0$. By continuity, there is some $\mu \in wL$ such that $B_\mu(\bar{c}d)$ contains no roots of $g_\alg$ for each of the $g$ appearing above, and also does not contain any roots of $s(F)_\alg$. Further, by Lemma \ref{valuation_atomics}, there is some $\lambda \in wL$ such that for any $\bar{t}u \in L$ with $\Jet(\bar{t})\Jet(u) \in B_\lambda(\bar{c}\bar{d})$, $\chi(\bar{t}u)$ holds. 

Let $M$ be the highest derivative of $x$ appearing in the formula $\phi$. Since the differential algebraic relation $F(\bar{y}, x) = 0$ holds, we may express the derivatives $y^{(m+1)},...,y^{(M)}$ as continuous (in fact, rational) functions of $\bar{y}$ and $x,x',...,x^{(m)}$. Therefore, there is some $\zeta \in wL$ such that, for any $\bar{t}u \in L$ with $F(\bar{t}, u) = 0$, $s(F)(\bar{t}, u) \neq 0$, $\Jet(\bar{t}) \in B_\zeta(\bar{b})$ and $\Jet_m(u) \in B_\zeta(\bar{c}\rest_m)$, we have that $\Jet_M(u) \in B_{\max(\lambda, \mu, \gamma)}(\bar{c}\rest_M)$. 

By the implicit function theorem for henselian fields (Theorem \ref{henselian_ift}), there is some $\theta \in wL$, such that for any differentially algebraically independent tuple $\bar{t}$ over $A$ with $\Jet(\bar{t}) \in B_\theta(\bar{c})$, there is a unique $\hat{d}_m \in L$ such that $F_\alg(\Jet(\bar{t}), d_0,...,d_{m-1},\hat{d}_m) = 0$, and $w(\hat{d}_m - d_m) > \max(\lambda, \mu, \zeta, \gamma) = \eta$. 

Let $\epsilon = \max(\lambda, \mu, \theta, \zeta, \gamma)$. By Lemma \ref{embedding_transcendentals}, there is a tuple $\bar{t} \in L$ differentially algebraically independent over $A$, realising $\qf\tp_\vd(\bar{a}/A)$, and $\Jet(\bar{t}) \in B_\epsilon(\bar{c})$. We find $\hat{d}_m \in B_\eta(d_m)$ such that 
\[
F_\alg(\Jet(\bar{t}), d_0,...,d_{m-1},\hat{d}_m) = 0
\]
and
\[
s(F)_\alg(\Jet(\bar{t}), d_0,...,d_{m-1},\hat{d}_m) \neq 0.
\]
Therefore, by differential henselianity, there is a $u \in L$ such that $F(\bar{t}, u) = 0$, and $\Jet_M(u) \in B_\eta(\bar{d}\rest_M)$.

Thus, $L \models \phi(\bar{t}, u)$, and the partial type $p$ is finitely satisfiable in $L$. By saturation, we find $\bar\alpha\beta \in L$ realising $p$, as required. 
\end{proof}

We will now show that the valued-differential quantifier-free type of a arbitrary finite tuple with a realisation of its valued-field quantifier-free type can be realised.

\begin{nota}
    Let $\L$ be an arbitrary language, and let $K, L$ be $\L$-structures with a common subset $A$. We write
    \[
    K \eqra_{\exists, A} L
    \]
    if every existential $\L(A)$-sentence which holds in $K$ also holds in $L$.
\end{nota}

\begin{prop} \label{approximating_types}
    Let $(K, v, \delta)$ be a valued-differential field, $(L, w, \partial) \models \DHen$, and $A$ a common valued-differential subfield. Suppose that $(L, w, \partial)$ is $|A|^+$-saturated and that $(K, v) \eqra_{\exists, A} (L, w)$ as pure valued fields. Let $\bar{a} = (a_i)_{i<n} \in K$ be an arbitrary finite tuple. Suppose there is $\bar{b} = (b_{ij})_{i<n, j<\omega}$ realising $\qf\tp_\vf(\Jet(\bar{a})/A)$. Then, for any $\gamma \in wL$, there is $\bar{c} \in L$ realising $\qf\tp_\vd(\bar{a}/A)$ such that $\Jet(\bar{c}) \in B_\gamma(\bar{b})$.
\end{prop}
\begin{proof}
    We will show that the quantifier-free type $p(\bar{x}) = \qf\tp_\vd(\bar{a}/A)$ along with the partial type stating that `$\Jet(\bar{x}) \in B_\gamma(\bar{b})$' is finitely satisfiable in $(L, w, \partial)$, and use compactness.

    First, we assume that the differential field extension $K \subseteq K\la \bar{a}\ra = F$ has nonzero differential transcendence degree. If not, possibly after replacing $K$ with a suitable elementary extension, adjoin an arbitrary differentially transcendental element $t$ to $\bar{a}$.

    Taking a differential transcendence basis, we can find $\bar{a}_0$ and $\bar{a}_1$ partitioning $\bar{a}$ such that $\bar{a}_0$ is differentially algebraically independent over $A$, and the extension $A\la \bar{a}_0 \ra \subseteq A\la \bar{a}_0 \ra \la \bar{a}_1 \ra = F$ is differentially algebraic. Partition the tuple $\bar{x}$ of variables similarly as $\bar{x}_0 \bar{x}_1$ and similarly reindex $\bar{b} = (b_{ij})$ as $\bar{b}_0 \bar{b}_1 = (b_{0, i, j})(b_{1, i, j})$. That is, $\bar{b}_0\bar{b}_1$ realises $\qf\tp_\vf(\Jet(\bar{a}_0)\Jet(\bar{a}_1)/A)$.

    As the derivative on $A\la \bar{a}_0\ra$ is nontrivial (in particular, there exists an element differentially transcendental over $A$), we may apply the differential primitive element theorem (\cite[Proposition II.9]{Kolchin1973}) to obtain a single element $c \in F$ such that $A \la \bar{a}_0 \ra\la c\ra = F$. 
    In particular, writing $\bar{a}_1 = (a_{1, i})_{i<m}$, each derivative $a_{1,i}^{(j)}$ is expressible as $q_{i,j}(\bar{a}_0, c)$, where $q(\bar{x}_0, y)$ is a differential rational function over $A$ (i.e. a ratio of differential polynomials over $A$ with non-vanishing denominator).

    As $\bar{b}$ is a realisation of $\qf\tp_\vf(\bar{a}/A)$, setting $\Jet(\bar{a}) \mapsto \bar{b}$ induces an embedding of valued fields over $A$. Let $\bar{d}$ be the image of $\Jet(c)$ under this embedding, i.e. $\bar{b}_0\bar{d}$ is a realisation of $\qf\tp_\vf(\Jet(\bar{a}_0c))$ in $L$.

    % % In particular, we can express each element $a_{1,i}$ of $\bar{a}_1 = (a_{1, i})_{i<m}$ as a rational combination over $A$ of finitely many derivatives of the members of $\bar{a}_0$ and $\alpha$. Further, every derivative of $a_{1, i}$ is also a rational combination of the derivatives of $\bar{a}_0$ and $\alpha$.
    % Let $\phi \in p$ be a formula. As $\phi$ is finite, only finitely many derivatives of the $\bar{x}_1 = (x_{1,i})_{i<m}$ appear in $\phi$.
    By continuity of the $q_{i,j}$, for any $N<\omega$, there is some $\delta \in wL$ such that for any realisation $\bar\alpha\beta$ of $\qf\tp_\vd(\bar{a}_0c/A)$ in $L$ with $\Jet(\bar\alpha\beta) \in B_\delta(\bar{b}_0\bar{d})$, we have that $w(q_{i, j}(\bar\alpha, \beta) - b_{1, i, j}) > \gamma$ for each $i<m, j<N$.

    By Lemma \ref{approximating_types_restricted_form}, there is a realisation $\bar\alpha\beta$ of $\qf\tp(\bar{a_0}c)$ in $L$ such that $\Jet(\bar\alpha \beta) \in B_{\max(\delta, \gamma)}(\bar{b}_0\bar{d})$. Let $\phi: A\la \bar{a}\ra \to L$ be the valued-differential field embedding induced by setting $\bar{a}_0c \mapsto \bar\alpha \beta$. By construction, $\phi(\Jet_N(\bar{a}_0\bar{a}_1)) \in B_\gamma(\bar{b}_0\bar{b}_1)$.

    Thus, $\phi(\bar{a}_0\bar{a}_1)$ is a realisation in $(L, w, \partial)$ of $\qf\tp_\vd(\bar{a}_0\bar{a}_1/A)$ as well as the formula stating `$\Jet_N(\bar{x}_0\bar{x}_1) \in B_\gamma(\bar{b}_0\rest_N \bar{b}_1\rest_N)$. We conclude therefore that the desired partial type is finitely satisfiable in $(L, w, \partial)$.
\end{proof}

From this, we harvest a number of more applicable results regarding embeddings of valued-differential fields.

\begin{prop}\label{existential_embedding}
% Let $(K, v, \delta)$ be a valued-differential field, $(L, w, \partial) \models \DHen$, and $A$ a common valued-differential subfield. Let $(L, w, \partial)$ be $|A|^+$-saturated. 
Suppose that $(K, v) \equiv\ra_{\exists, A} (L, w)$ as pure valued fields. Then, every differentially finitely generated extension $A\la \bar{a}\ra \subseteq K$ embeds in $L$ over $A$ as valued-differential fields.
\end{prop}
\begin{proof}
% We may assume that the tuple $\bar{a}$ is of the form $\bar{a}_0\bar{b}$, where $\bar{a}_0$ is differentially algebraically independent over $A$, and the extension $A\la \bar{a}_0\bar{b}\ra$ is differentially algebraic over $A\la\bar{a}_0\ra$.

% We may also assume that $A\la a_0\ra$ is not a constant field, i.e. $A$ is not constant or $a_0$ is a nonempty tuple. Otherwise, we take an element $t \in K$ (possibly after replacing $K$ with an elementary extension) differentially transcendental over $A\la \bar{b}\ra$, and considering the extension $A\la \bar{a}t\ra$. 

% Since $K \equiv\ra_{\exists, A} L$ as valued fields, and by saturation of $L$, there is a realisation $\bar{c}\bar{d} \in L$ of $\qf\tp_\vf(\Jet(\bar{a_0})\Jet(b)/A)$. Now, by Lemma \ref{approximating_types}, there is a realisation $\bar\alpha\beta$ of $\qf\tp_\vd(\bar{a}_0b/A)$. Setting $\bar{a}_0b \mapsto \bar\alpha\beta$, we obtain an embedding of $A\la \bar{a}\ra \to L$ as valued-differential fields.

As $(K, v) \eqra_{\exists, A} (L, w)$, and by saturation of $(L, w)$, there is a realisation $\bar{b}$ of $\qf\tp_\vf(\Jet(\bar{a}/A)$ in $L$. By Proposition \ref{approximating_types}, there is a realisation $\bar{c}$ of $\qf\tp_\vd(\bar{a}/A)$ in $(L, w, \partial)$. In particular, setting $\bar{a} \mapsto \bar{c}$ induces an embedding of valued-differential fields $A \la \bar{a} \ra \to L$, as required.
\end{proof}

From the above embedding lemma, we obtain the desired existential lifting property of differentially henselian fields (cf. \cite[Theorem 3.14]{Guzy2006}). Note that in the following, there is no saturation requirement on $(L, w, \partial)$.

\begin{cor}\label{eqra_valfield_implies_eqra_valdiff}
Let $(K, v, \delta)$ be a valued-differential field, $(L, w, \partial) \models \DHen$, and $A$ a common valued-differential subfield. Suppose that $(K, v) \equiv\ra_{\exists, A} (L, w)$ as pure valued fields. Then, $(K, v, \delta) \equiv\ra_{\exists, A} (L, w, \partial)$ as valued-differential fields.
\end{cor}
\begin{proof}
Let $\phi(\bar{x})$ be an $\Lvd(A)$-formula, and suppose that $(K, v, \delta) \models \exists\bar{x}\phi(\bar{x})$. Let $\bar{a} \in K$ such that $K \models \phi(\bar{a})$. Let $(L^+, w^+, \partial^+)$ be an $|A|^+$-saturated elementary extension of $(L, w, \partial)$. By Lemma \ref{existential_embedding}, there is an embedding $f: A\la \bar{a} \ra \to L^+$ as valued-differential fields. Since $\phi$ is quantifier-free and $f$ is an embedding of valued-differential fields, we have that $(L^+, w^+, \partial^+) \models \phi(f(\bar{a}))$, and so $(L^+, w^+, \partial^+) \models \exists\bar{x}\phi(\bar{x})$. Since $(L, w, \partial) \preceq (L^+, w^+, \partial^+)$, we also have that $(L, w, \partial) \models \exists\bar{x} \phi(\bar{x})$ as required.
\end{proof}

\begin{cor} 
Let $(K, v, \delta)$ be a valued-differential field, $(L, w, \partial) \models \DHen$, and $A$ a common valued-differential subfield. Suppose that $(K, v) \equiv\ra_{\exists,A} (L, w)$ as pure valued fields, and $(L, w, \partial)$ is $\max(|A|^+, |K|)$-saturated. Then, $(K, v, \delta)$ embeds in $(L, w, \partial)$ over $A$ as valued differential fields. 
\end{cor}
\begin{proof}
Let $|K| = \kappa$, and enumerate $K$ as $\bar{a} = (a_\alpha)_{\alpha<\kappa}$. It suffices to find a realisation of the quantifier-free type $p(\bar{x}) = \qf\tp_\vd(\bar{a}/A)$. As $(K, v) \eqra_{\exists, A} (L, w)$ as pure valued fields and $(L, w, \partial)$ is differentially henselian, by Corollary \ref{eqra_valfield_implies_eqra_valdiff}, we also have that $(K, v, \delta) \eqra_{\exists, A} (L, w, \partial)$ as valued-differential fields. Thus, $p$ is finitely satisfiable in $L$, and by saturation, we also have that $p$ is realised in $L$.
\end{proof}

% \begin{cor} \label{diff_henselian_ec_val_implies_ec_diffval}
% Let $(K, v, \delta)$ be a differentially henselian field, and let $(L, w, \partial)$ be a valued-differential field extension of $(K, v, \delta)$. If $(K, v)$ is existentially closed in $(L, w)$ as a valued field, then $(K, v, \delta)$ is also existentially closed in $(L, w, \partial)$ as a valued-differential field.
% \end{cor}
% \begin{proof}
% Observe that $K$ being existentially closed in $L$ as a valued field is equivalent to having $L \equiv\ra_{\exists, K} K$ as valued fields. By Corollary \ref{eqra_valfield_implies_eqra_valdiff}, we have that $L \equiv\ra_{\exists, K} K$ as valued-differential fields. Thus $K$ is existentially closed in $L$ as valued-differential fields.
% \end{proof}

We now prove the reverse direction of Theorem \ref{dh_equiv_model_DH_thm}, i.e. the converse of Lemma \ref{dh_implies_model_DH_lem}.

\begin{lem} \label{model_DH_implies_dh_lem}
Let $(K, v, \delta) \models \DHen$ be a valued-differential field. Then $(K, v, \delta)$ is differentially henselian.
\end{lem}
\begin{proof}
Let $(L, w, \partial)$ be a valued-differential field extension such that $(K, v)$ is existentially closed in $(L, w)$ as pure valued fields. Observe that $(K, v) \preceq_\exists (L, w)$ is equivalent to the condition $(L, w) \eqra_{\exists, K} (K, v)$. By Corollary \ref{eqra_valfield_implies_eqra_valdiff}, we also have that $(L, w, \partial) \eqra_{\exists, K} (K, v, \delta)$ as valued-differential fields, i.e. $(K, v, \delta)$ is existentially closed in $(L, w, \partial)$.
\end{proof}

This establishes Theorem \ref{dh_equiv_model_DH_thm}. In fact, it is possible to refine the statement of Theorem \ref{dh_equiv_model_DH_thm} to reference only irreducible differential polynomials.

\begin{lem} \label{factorisation_lem_1}
    Let $(K, \delta)$ be a differential field, and $f(x) \in K\{x\}$ be a differential polynomial of order $n$ such that $f = gh$, where $g, h \in K\{x\}$ are not units. Let $\bar{a} \in K^{n+1}$ such that $f_\alg(\bar{a}) = 0$ and $s(f)_\alg(\bar{a}) \neq 0$. Then, either
    \begin{enumerate}[label=(\roman*)]
        \item $g_\alg(\bar{a}) = 0, h_\alg(\bar{a}) \neq 0, s(g)_\alg(\bar{a}) \neq 0$ and $\ord(g) = n$, or
        \item $h_\alg(\bar{a}) = 0, g_\alg(\bar{a}) \neq 0, s(h)_\alg(\bar{a}) \neq 0$ and $\ord(h) = n$.
    \end{enumerate}
\end{lem}
\begin{proof}
    As $f_\alg(\bar{a}) = 0$ and $f = gh$, at least one of $g_\alg(\bar{a})$ and $h_\alg(\bar{a})$ vanishes. Consider $s(f)_\alg(\bar{a}) = \frac{\partial f}{\partial x^{(n)}}(\bar{a})$. By the product rule,
    \[
    s(f)_\alg(\bar{a}) = \left(\frac{\partial g_\alg}{\partial x^{(n)}}\right)(\bar{a}) h_\alg(\bar{a}) + g_\alg(\bar{a}) \left(\frac{\partial h_\alg}{\partial x^{(n)}}\right)(\bar{a}) \neq 0
    \]
    In particular, it is not the case that both $g_\alg(\bar{a})$ and $h_\alg(\bar{a})$ vanish, otherwise $s(f)_\alg(\bar{a}) = 0$, a contradiction. Without loss, suppose that $g_\alg(\bar{a}) = 0$ and $h_\alg(\bar{a}) \neq 0$. If $\ord(g) < n$, then $\frac{\partial g_\alg}{\partial x^{(n)}} = 0$, and $s(f)_\alg(\bar{a}) = 0$, again a contradiction. Thus $\ord(g) = n$. Now,
    \[
    s(f)_\alg(\bar{a}) = s(g)_\alg(\bar{a}) h_\alg(\bar{a}) + 0 \neq 0
    \]
    and since $h_\alg(\bar{a}) \neq 0$, we obtain that $s(g)_\alg(\bar{a}) \neq 0$ also. The other case follows by symmetry.
\end{proof}

\begin{lem} \label{factorisation_lem_2}
    Let $(K, \delta)$ be a differential field, and $f(x) \in K\{x\}$ be a differential polynomial of order $n$, and $\bar{a} \in K^{n+1}$ such that $f_\alg(\bar{a}) = 0$ and $s(f)(\bar{a}) \neq 0$. Then, there is a unique irreducible factor $g$ of $f$ such that $g_\alg(\bar{a}) = 0$. Further, $\ord(g) = n$ and $s(g)_\alg(\bar{a}) \neq 0$.
\end{lem}
\begin{proof}
    As $K\{x\}$ is a unique factorisation domain, we may write $f$ as the product of its irreducible factors. Applying Lemma \ref{factorisation_lem_1} repeatedly to this factorisation, we observe that there is a unique irreducible factor $g$ such that $g$ vanishes on $\bar{a}$. Further, also by Lemma \ref{factorisation_lem_1}, such a factor necessarily has order $n$, and $s(g)_\alg(\bar{a}) \neq 0$.
\end{proof}

Now, it is clear that we can restrict our axiomatisation of differentially henselian fields to only irreducible polynomials:

\begin{prop} \label{dh_iff_dh_for_irreducible_prop}
    Let $(K, v, \delta)$ be a valued-differential field, henselian as a pure valued field. Then, $(K, v, \delta)$ is differentially henselian if and only if it satisfies the following axiom scheme:
    \begin{center}
        For every irreducible differential polynomial $f(x) \in K\{x\}$ of order $n$, and $\bar{a} \in K^{n+1}$ such that $f_\alg(\bar{a}) = 0$ and $s(f)_\alg(\bar{a}) \neq 0$, for any $\gamma \in vK$, there is $b \in K$ such that $f(b) = 0$ and $\Jet_n(b) \in B_\gamma(\bar{a})$.
    \end{center}
\end{prop}
\begin{proof}
    We prove the condition established in Theorem \ref{dh_equiv_model_DH_thm}. The forwards direction is trivial. For the reverse, let $f(x) \in K\{x\}$ be an arbitrary differential polynomial of order $n$ and $\bar{a} \in K^{n+1}$ such that $f_\alg(\bar{a}) = 0$ and $s(f)_\alg(\bar{a}) \neq 0$. Let $\gamma \in vK$.

    By Lemma \ref{factorisation_lem_2}, there is a unique irreducible factor $g$ of $f$ such that $g_\alg(\bar{a}) = 0$, $s(g)_\alg(\bar{a}) \neq 0$ and $\ord(g) = n$. By the axiom scheme above, there is $b \in K$ such that $g(b) = 0$ and $\Jet_n(b) \in B_\gamma(\bar{a})$. Then as $g$ is a factor of $f$, we also have that $f(b) = 0$, as required.
\end{proof}

We now harvest a useful relative embedding theorem for differentially henselian fields.

\begin{thm}[Relative Embedding Theorem]\label{relative_embedding_thm} 
Let $(K, v, \delta)$ and $(L, w, \partial)$ be valued-differential fields, where $(L, w, \partial)$ is differentially henselian, and let $A$ a common valued-differential subfield. Suppose that there is an embedding $\phi: (K, v) \to (L, w)$ as pure valued fields over the common subfield $A$. Suppose also that $(L, w, \partial)$ is $|K|^+$-saturated. Then, for any $\gamma \in wL$, there is an embedding $\psi: (K, v, \delta) \to (L, w, \partial)$ of valued-differential fields over $A$ such that for all $a \in K$, $\psi(a) \in B_\gamma(\phi(a))$. 
\end{thm}
\begin{proof}
Enumerate $K$ as $\bar{k}$ and consider the quantifier-free type $p(\bar{x})$ which contains $\qf\tp_\vd(\bar{k}/A)$ along with the formulae which express `$\bar{x} \in B_\gamma(\phi(\bar{k})$'. This is finitely satisfiable by Proposition \ref{approximating_types}, as any finite subset of $p$ involves only finitely many variables, and $p$ restricted to these variables is realisable in $L$. Applying compactness, we find a realisation $\bar{l}$ of $p$ in $L$, i.e. setting $\bar{k} \mapsto \bar{l}$ induces
a differential embedding $\psi$ of $K$ into $L$ such that for any $a \in K$, we have that $\psi(c) \in B_\gamma(\phi(a))$.
\end{proof}

We will now demonstrate that the class of differentially henselian fields coincides with the henselian models of Guzy's $(\UC'_k)$ when $k=1$. The argument is essentially the same as for the differentially large case, i.e. Proposition 4.7 of \cite{SanchezTressl2020}.

\begin{prop}
    Let $(K, v, \delta)$ be a valued-differential field, henselian as a pure valued field. Then, $(K, v, \delta)$ is differentially henselian if and only if $(K, v, \delta) \models (\UC'_1)$. 
\end{prop}
\begin{proof}
    First suppose that $(K, v, \delta)$ is differentially henselian. Then, by \cite[Theorem 3.14]{Guzy2006}, there is a valued-differential field extension $(L, w, \partial)$ of $(K, v, \delta)$ such that $(L, w, \partial) \models (\UC'_1)$ and $(K, v) \preceq (L, w)$ as pure valued fields. Then, as $(\UC'_1)$ is inductive, $(K, v, \delta)$ is also a model of $(\UC'_1)$.

    Now, suppose that $(K, v, \delta) \models (\UC'_1)$. Let $(L, w, \partial)$ be any valued-differential field extension of $(K, v, \delta)$ such that $(K, v) \preceq_\exists (L, w)$. By existential closure, there is an elementary extension $(F, u)$ of $(K, v)$ such that $(L, w)$ embeds in $(F, u)$ over $(K, v)$. Extend the derivation $\partial$ arbitrarily to a derivation $\dd$ on $F$. It suffices to show that $(K, v, \delta)$ is existentially closed in $(F, u, \dd)$, so replace $(L, w, \partial)$ with $(F, u, \dd)$ (in particular, we may assume that $(K, v) \preceq (L, w)$).

    Again by Theorem \cite[Theorem 3.14]{Guzy2006}, there is a valued-differential field extension $(L^*, w^*, \partial^*)$ of $(L, w, \partial)$ such that $(L^*, w^*, \partial^*) \models (\UC'_1)$ and $(L, w) \preceq (L^*, w^*)$. In particular, $(K, v) \preceq (L^*, w^*)$ as valued fields, and thus $(K, v) \preceq_\exists (L^*, w^*)$, i.e. $(K, v)$ and $(L, w)$ have the same universal theory over $K$ as pure valued fields. 
    Now, as both $(K, v, \delta)$ and $(L^*, w^*, \partial^*)$ are models of $(\UC'_1)$, and they have the same universal $\Lvf$-theory over $K$, by \cite[Theorem 3.14]{Guzy2006}, they have the same universal $\Lvd$-theory over $K$. That is, $(K, v, \delta)$ is existentially closed in $(L^*, w^*, \partial^*)$. This is sufficient, as $(K, v, \delta) \subseteq (L, w, \partial) \subseteq (L^*, w^*, \partial^*)$, so $(K, v, \delta) \preceq_\exists (L, w, \partial)$, as required. 
\end{proof}

\section{Ax-Kochen/Ershov Principles for Differentially Henselian Fields} \label{section_dh_ake}

In this section, we will apply the embedding theorem from the previous section to show that many Ax-Kochen/Ershov type results can be lifted from henselian fields to differentially henselian fields.

\begin{defn}
    Let $\mathcal{C}$ be an elementary class of henselian valued fields or differentially henselian fields. We say that $\mathcal{C}$ is
    \begin{itemize}
        \item \textbf{relatively complete}, or a \textbf{$\AKE^\equiv$-class}, if for any $(K, v), (L, w) \in \mathcal{C}$ if $vK \equiv wL$ and $Kv \equiv Lw$, then $(K, v) \equiv (L, w)$;
        % \item \textbf{relatively subcomplete} if for any $(K, v), (L, w) \in \mathcal{C}$ with a common valued subfield $(F, u)$, if $vK \equiv_{uF} wL$ and $Kv \equiv_{Fu} Lw$, then $(K, v) \equiv_F (L, w)$;
        \item \textbf{relatively model complete}, or a \textbf{$\AKE^\prec$-class} if for any $(K, v) \subseteq (L, w) \in \mathcal{C}$, if $vK \preceq wL$ and $Kv \preceq Lw$, then $(K, v) \preceq (L, w)$;
        \item \textbf{relatively existentially complete}, or \textbf{a $\AKE^\exists$-class} if for any $(K, v) \subseteq (L, w) \in \mathcal{C}$, if $vK \preceq_\exists wL$ and $Kv \preceq_\exists Lw$, then $(K, v) \preceq_\exists (L, w)$.
    \end{itemize}
    For a class $\mathcal{C}$ of differentially henselian fields, we denote the class of henselian valued fields obtained by taking reducts of members of $\mathcal{C}$ by $\mathcal{C}_\alg$.
\end{defn}

\begin{egs}
    \begin{enumerate}
        \item The class of henselian valued fields of equicharacteristic 0 is relatively complete, model complete, and existentially complete.
        \item The class of algebraically closed valued fields (of any characteristic) is relatively complete, model complete, and existentially complete.
        \item The class of unramified henselian valued fields of mixed characteristic $(0, p)$ (i.e. $(K, v)$ is unramified if $vp$ is minimal positive in $vK$) is relatively complete and model complete \cite[Corollary 8.5, Theorem 9.2]{AnscombeJahnke2019}, and is relatively existentially complete under the additional condition that every residue field has the same finite degree of imperfection \cite[Theorem 10.2]{AnscombeJahnke2019}.

    \end{enumerate}
\end{egs}

For the remainder of this section, we let $\mathcal{C}$ be a class of differentially henselian fields.

\begin{prop} \label{ake_exists_prop}
    If $\mathcal{C}_\alg$ is an $\AKE^\exists$-class, then $\mathcal{C}$ is also an $\AKE^\exists$-class.
\end{prop}
\begin{proof}
    This follows immediately from the definition. Let $(K, v, \delta) \subseteq (L, w, \partial) \in \mathcal{C}$. Suppose that $vK \preceq_\exists wL$ and $Kv \preceq_\exists Lw$. As $\mathcal{C}_\alg$ is a $\AKE^\exists$-class, $(K, v) \preceq_\exists (L, w)$. Since $(K, v, \delta)$ is differentially henselian, $(K, v, \delta) \preceq_\exists (L, w, \partial)$.
\end{proof}

\begin{prop} \label{ake_equiv_prop}
    If $\mathcal{C}_\alg$ is an $\AKE^\equiv$-class, then $\mathcal{C}$ is also an $\AKE^\equiv$-class.
\end{prop}
\begin{proof}
    This follows from \cite[Corollary 2.4.7]{CubidesPoint2019}. Let $(K, v, \delta),(L, w, \partial) \in \mathcal{C}$, and suppose that $vK \equiv wL$ and $Kv \equiv Lw$. As $\mathcal{C}_\alg$ is an $\AKE^\equiv$-class, we have that $(K, v) \equiv (L, w)$. By \cite[Corollary 2.4.7]{CubidesPoint2019}, the theory of a differentially henselian field whose underlying valued field is elementarily equivalent to $(K, v)$ is complete. Thus $(K, v, \delta) \equiv (L, w, \partial)$.
\end{proof}

\begin{thm} \label{ake_prec_thm}
    If $\mathcal{C}_\alg$ is an $\AKE^\prec$-class, then $\mathcal{C}$ is also an $\AKE^\prec$-class.
\end{thm}
\begin{proof}
    We will prove that for any two differentially henselian fields $(K, v, \delta)$ and $(L, w, \partial)$ in $\mathcal{C}$ with a common valued-differential subfield $(F, u, \dd) \in \mathcal{C}$ satisfying $Fu \preceq Kv, Lw$ and $uF \preceq vK, wL$, we have that $(K, v, \delta) \equiv_F (L, w, \partial)$. By taking appropriate elementary extensions, we may assume that $(K, v, \delta)$ and $(L, w, \partial)$ are $\kappa=|F|^+$-saturated.

    Let $(K_0, v_0, \delta_0)$ be an elementary substructure of $(K, v, \delta)$ over $F$ of cardinality $\kappa$. As $\mathcal{C}_\alg$ is an $\AKE^\prec$-class, we have that $(K_0, v_0)_F \equiv (L, w)_F$. By saturation, there is an embedding $\phi_0: (K_0, v_0) \to (L, w)$ over $F$. Applying Theorem \ref{relative_embedding_thm}, there is a differential embedding $\psi_0: (K_0, v_0, \delta_0) \to (L, w, \partial)$ over $F$ such that $\psi_0$ and $\phi_0$ induce the same embeddings of the value group and residue field. Denote by $(L_0, w_0, \partial_0)$ the image of $\psi_0$. Note that $\psi_0$ remains a elementary embedding of valued fields. We inductively construct the isomorphisms $\psi_n: (K_n, v_n, \delta_n) \to (L_n, w_n, \partial_n)$ for all $n$.

    Suppose that we have constructed subfields $(K_{2n}, v_{2n}, \delta_{2n})$ and $(L_{2n}, w_{2n}, \partial_{2n})$ of $(K, v, \delta)$ and $(L, w, \partial)$, respectively of cardinality $\kappa$, and an isomorphism $\psi_{2n}: (K_{2n}, v_{2n}, \delta_{2n}) \to (L_{2n}, w_{2n}, \partial_{2n})$ such that $(K_{2n}, v_{2n}, \delta_{2n}) \preceq (K, v, \delta)$, and $\psi_{2n}$ is an elementary embedding of valued fields.

    Let $(L_{2n+1}, w_{2n+1}, \partial_{2n+1})$ be a elementary substructure of $(L, w, \partial)$ of cardinality $\kappa$ over $L_{2n}$. Since $(L_{2n+1}, w_{2n+1})_{L_{2n}} \equiv (K, v)_{K_2n}$, by saturation, there is an elementary embedding $\phi^{-1}_{2n}$ of $(L_{2n+1}, w_{2n+1})$ into $(K, v)$ over $K_{2n}$. By Theorem \ref{relative_embedding_thm}, there is an embedding $\psi_{2n+1}^{-1}: (L_{2n+1}, w_{2n+1}, \partial_{2n+1}) \to (K, v, \delta)$ over $K_{2n}$, which induces the same embeddings of value groups and residue fields as $\phi_{2n+1}^{-1}$. Note that this implies $\psi_{2n+1}^{-1}$ is an elementary embedding of valued fields. Let $(K_{2n+1}, v_{2n+1}, \delta_{2n+1})$ be the image of $\psi_{2n+1}^{-1}$.

    Exchanging the roles of $L$ and $K$, similarly construct the isomorphism $\psi_{2n+2}: (K_{2n+2}, v_{2n+2}, \delta_{2n+2}) \to (L_{2n+2}, w_{2n+2}, \partial_{2n+2})$. Let
    \begin{align}
        (K_\infty, v_\infty, \delta_\infty) = \bigcup_{n<\omega} (K_n, v_n, \delta_n), \\
        (L_\infty, w_\infty, \partial_\infty) = \bigcup_{n<\omega} (L_n, w_n, \partial_n).
    \end{align}
    Observe that $K_\infty$ and $L_\infty$ are isomorphic as valued-differential fields over $F$, and that
    \[
    (K, v, \delta)_F \equiv (K_\infty, v_\infty, \delta_\infty)_F \equiv (L_\infty, w_\infty, \partial_\infty)_F \equiv (L, w, \partial)_F
    \]
    as required.
\end{proof}

\begin{defn}
    Following Kuhlmann, we say the class $\mathcal{C}$ is \textbf{relatively subcomplete} if for every pair $(K, v), (L, w) \in \mathcal{C}$ and common valued(-differential) subfield $(F, u)$ such that
    \begin{enumerate}[label=(\roman*)]
        \item $(F, u)$ is defectless;
        \item $vK/uF$ is torsion free, i.e. $uF$ is \textbf{pure} in $vK$;
        \item $Kv/Fu$ is separable;
    \end{enumerate}
    if
    \[
    vK \equiv_{uF} wL \text{ and } Kv \equiv_{Fu} Lw
    \]
    then $(K, v) \equiv_F (L, w)$.
\end{defn}
We direct the reader to \cite{Kuhlmann2016} for the definition of the defect.

\begin{prop}
    Let $\mathcal{C}$ be an elementary class of defectless differentially henselian fields such that $\mathcal{C}_\alg$ is relatively subcomplete. Then, $\mathcal{C}$ is relatively subcomplete.
\end{prop}
\begin{proof}
    This is essentially the same as the proof of Theorem \ref{ake_prec_thm}, making appropriate adjustments.
\end{proof}

\section{Relative Quantifier Elimination in the Pas Language} \label{section_dh_qe_for_ac}

In this section, we adapt relative quantifier elimination results for equicharacteristic 0 henselian valued fields in the Pas language to the differential context. We generalise results by Borrata (\cite[Corollary 4.3.27]{Borrata2021}) on closed ordered differential valued fields to arbitrary differentially henselian fields. We begin by recalling a definition:

\begin{defn}
Let $(K, v)$ be a valued field. An \textbf{angular component map} for $(K, v)$ is a map $\ac: K \to Kv$ satisfying:
\begin{enumerate}[label=(\arabic*)]
\item $\ac(x) = 0$ if and only if $x = 0$;
\item $\ac\rest_{K^\times}: K^\times \to Kv^\times$ is a group homomorphism, and
\item $\ac\rest_{\O_v^\times} = \res\rest_{\O_v^\times}$.
\end{enumerate}
An $\ac$-valued field is a valued field equipped with an angular component map. Similarly, an $\ac$-valued differential field is an $\ac$-valued field equipped with a derivation.
\end{defn}

\begin{eg}
Let $K$ be a field of characteristic 0, and consider the field of formal Laurent series $K((t))$. For a nonzero $a \in K((t))$, we set $\ac(a) = a_{v(a)}$, i.e. $\ac(a)$ is the first non-zero coefficient in $a$. 
\end{eg}

\begin{note}
Not every henselian valued field admits an angular component map.
\end{note}

If we have a cross section $s: vK \to K^\times$, that is, $s$ is a group homomorphism satisfying $v\circ s = \id_{vK}$, then we can define an angular component map by setting
\[
\ac(x) = \left\{ \begin{array}{cc}
\res(x/s(v(x))) & \text{for } x \neq 0 \\
0 & \text{otherwise}.
\end{array}\right.
\]

\begin{eg}
The angular component map for $K((t))$ in the example above is given by the cross-section $s: \Z \to K((t))$ defined by $n \mapsto t^n$.
\end{eg}

\begin{lem}
Let $(K, v) \subseteq (L, w)$ be an unramified extension of valued fields, i.e. $vK = wL$, and let $\ac$ be an angular component map for $(K, v)$. Then, there is a unique angular component map on $L$ extending $\ac$.
\end{lem}
\begin{proof}
Let $b \in L^\times$. Since the extension is unramified, there is $a \in K^\times$ such that $v(a) = w(b)$. Then, $w(b/a) = 0$, and thus $b/a \in \O_w^\times$. If $\ac$ is an angular component map on $(L, w)$, then we have that $\ac(b/a) = \res(b/a)$. By multiplicativity of $\ac$, we have that $\ac(b) = \ac(b/a)\ac(a) = \res(b/a)\ac(a)$, which is uniquely determined by $\ac$ on $K$. This gives the unique extension of $\ac$ to $L$.
\end{proof}

We work in the three-sorted \textbf{Denef-Pas language} $\L_\Pas$, with a sort $K$ for the valued field (the `field sort'), a sort $\Gamma$ for the value group, and a sort $k$ for the residue field. We add the usual symbols for operations and relations on each of the sorts, i.e. the language of rings in the sorts $K$ and $k$, and the language of ordered abelian groups in the sort $\Gamma$. We also include symbols $v:K \to \Gamma$ for the valuation, and a symbol $\ac: K \to k$ for an angular component map. Note that the residue map is quantifier-free definable from $\ac$ as $\ac = \res$ on $\mathcal{O}^\times$, and we can set $\res(x) = 0$ otherwise.

In \cite{Pas1989}, Pas shows that the theory of henselian $\ac$-valued fields of equicharacteristic 0 eliminates $K$-quantifiers. We show that an analogous result holds in the valued-differential setting.

We denote by $\L_\Pas^\delta$ the language $\L_\Pas$ augmented with a unary function symbol $\delta: K \to K$ in the field sort. That is, $\L_\Pas^\delta$ is the \textbf{language of differential $\ac$-valued fields}.

From the proof of the classical theorem found in \cite{HilsNotes}, we can extract the following embedding lemma for henselian $\ac$-valued fields:
\begin{lem}\label{classical_ac_embedding}
Let $(K, v, \ac_K)$ and $(L, w, \ac_L)$ be henselian $\ac$-valued field, considered as $\L_\Pas$-structures. Suppose that $K$ is countable, and $(L, w, \ac)$ is $\aleph_1$-saturated. Let $(A, u, \ac_A)$ be an $\L_\Pas$-substructure of $(K, v, \ac_K)$, and let $f: (A, u, \ac_A) \to (L, w, \ac_L)$ be an $\L_\Pas$-embedding such that the induced maps $f_g$ and $f_r$ on the value groups and residue fields are elementary. Then, $f$ extends to an $\L_\Pas$-embedding of $(K, v, \ac_K)$ into $(L, w, \ac_L)$.
\end{lem}

We require one more preparatory lemma before we proceed to the relative quantifier elimination result:

\begin{lem} \label{ac_embedding_after_small_change_remains_ac_embedding_lem}
    Let $(K, v, \ac_K)$ and $(L, w, \ac_L)$ be $\ac$-valued fields. Let $\phi: K \to L$ be an $\ac$-valued field embedding, and $\psi: K \to L$ be an embedding of valued fields such that $w(\phi(a) - \psi(a)) > w(\phi(K))$ for any $a \in K$. Then, for any $b \in K$, $\ac_L(\phi(b)) = \ac_L(\psi(b))$, and $\psi$ is also an embedding of $\ac$-valued fields.
\end{lem}
\begin{proof}
    We begin by showing that, for $a \in L^\times$ and $b \in L$ such that $w(b) > w(a)$, we have that $\ac_L(a) = \ac_L(a+b)$. By multiplicativity, we have that:
    \[
    \ac_L(a+b)(\ac_L(a^{-1})) = \ac_L(1 + a^{-1}b)
    \]
    Since $w(a) < w(b)$, we have $w(a^{-1}b) = w(b) - w(a) > 0$, thus $w(1 + a^{-1}b) = 0)$ and $1+ a^{-1}b \in \O_w^\times$. Thus,
    \[
    \ac_L(1 + a^{-1}b) = \res_w(1+a^{-1}b) = 1
    \]
    and so $\ac_L(a+b) = \ac_L(a) = \ac_L(a^{-1})^{-1}$ as required.

    Now, let $\phi, \psi$ be as above. Observe that for any $b \in K$, $w(\phi(b) - \psi(b)) > w(\phi(b))$, thus by the above claim, $\ac_L(\psi(b)) = \ac_L(\phi(b))$. Since $\phi$ is an $\ac$-valued field embedding, so is $\psi$.
\end{proof}

We will now apply the above lemmas to adapt the proof of the classical Pas' theorem found in \cite{HilsNotes} to the differentially henselian case.

\begin{thm} \label{dh_ac_has_qe}
The $\L_\Pas^\delta$-theory of differentially henselian $\ac$-valued fields eliminates quantifiers in the field sort.
\end{thm}
\begin{proof}
Let $(K, v, \ac_K, \delta), (L, w, \ac_L, \partial)$ be differentially henselian $\ac$-valued fields. We assume $K$ is countable and $(L, w, \partial, \ac_L)$ is $\aleph_1$-saturated. Let $(A, u, \ac_A, \dd)$ be an $\L_\Pas^\delta$-substructure of $K$, and $f: A \to L$ be an $\L_\Pas^\delta$-embedding, such that the embeddings $f_r$ and $f_g$ of the residue field and value group sorts are elementary with respect to the languages of rings and ordered abelian groups, respectively. It now suffices to show that $f$ extends to an $\L_\Pas^\delta$-embedding $K \to L$.

Now, by Lemma \ref{classical_ac_embedding}, there is an $\L_\Pas$-embedding $\phi: K \to L$ extending $f$. By $\aleph_1$-saturation of $L$, there is some $\gamma \in wL$ with $\gamma > w(\phi(K))$. Now, apply Theorem \ref{relative_embedding_thm} to obtain a valued-differential field embedding $\psi: K \to L$ extending $f$ with $w(\phi(a) - \psi(a)) > \gamma$ for every $a \in K$. By Lemma \ref{ac_embedding_after_small_change_remains_ac_embedding_lem}, $\psi$ also preserves $\ac$, thus is an $\L_\Pas^\delta$-embedding, as required.
\end{proof}

\section{Stable Embeddedness} \label{section_dh_stable_embedded}

In this section, we prove a stable-embeddedness result for the value group and residue field of a differentially henselian field in both the equicharacteristic 0 and unramified mixed characteristic cases, adapting analogous results for their pure valued field counterparts.

We recall the definition of stable-embeddedness:

\begin{defn}
Let $M$ be an arbitrary first-order structure, and let $P \subseteq M^n$ be a definable set. We say that $P$ is \emph{stably embedded} if, for all formulae $\phi(\bar{x}, \bar{y})$. and all $\bar{b} \in M^{|\bar{y}|}$, $\phi(M^{n|\bar{x}|}, \bar{b}) \cap P^{|\bar{x}|}$ is $P$-definable.
\end{defn}

To work with the value group and residue field as definable sets in the valued field, for the rest of this section, we consider a valued field $K$ as a three-sorted structure $(K, vK, Kv)$ equipped with the usual languages for each sort, along with functions for the valuation and residue map. Our result and proof are an adaptation of \cite[Theorem 11.3]{AnscombeJahnke2019} to the differential context, and by transitivity, also of \cite[Lemma 3.1]{JahnkeSimon2020}

\begin{thm} \label{dh_res_field_val_gp_stabl_emb_thm}
Let $(K, v, \delta)$ be an differentially henselian field, either of equicharacteristic 0 or mixed characteristic and unramified. Then, the value group $vK$ and residue field $Kv$ are stably embedded as a pure ordered abelian group and pure field, respectively.
\end{thm}
\begin{proof}
Let $(L, w, \partial)$ be a $|K|^+$-saturated elementary extension of $(K, v, \delta)$, and let $a, b \in Lw$ have the same type over $Kv$ in the language of rings. We will show that $a$ and $b$ have the same type over $K$.

Let $(L_0, w_0, \partial_0)$ be an $\aleph_1$-saturated elementary substructure of $(L, w, \partial)$ containing $Ka$. Since $a$ and $b$ have the same type over $Kv$, there is an elementary embedding $\phi_r: L_0w_0 \to L_w$ over $Kv$ with $\phi_r(a) = b$. Let $\phi_g: w_0L \to wL$ be the inclusion.

In the equicharacteristic 0 case, we observe that $vK$ is pure in $w_0L_0$, as $vK$ is elementary in $w_0L_0$. Thus we may apply Lemma \cite[Lemma 4.6.2]{PrestelDelzell2011} along with Theorem \ref{relative_embedding_thm} to find an embedding $\phi: (L_0, w_0, \partial_0) \to (L, w, \partial)$ inducing $\phi_r$ and $\phi_g$.

In the unramified mixed characteristic case, we apply observe that since $\phi_r$ is elementary, $Lw/\phi_r(L_0w_0)$ is a separable extension. Similarly to the equicharacteristic 0 case, $vK$ is pure in $w_0L_0$. Apply \cite[Proposition 10.1]{AnscombeJahnke2019} and Theorem \ref{relative_embedding_thm} to obtain an embedding $\phi: (L_0, w_0, \partial_0) \to (L, w, \partial)$ inducing $\phi_r$ and $\phi_g$.

Since the classes of differentially henselian fields of equicharacteristic 0 and unramified mixed characteristic are $\AKE_\prec$-classes, as $\phi_r$ and $\phi_g$ are elementary embeddings, so is $\phi$. 

We finally conclude the following:
\[
\tp_L(a/K) = \tp_{L_0}(a/K) = \tp_{\phi(L_0)}(b/K) = \tp_L(b/K),
\]
where $\tp_L(a/K)$ denotes the valued-differential field type of $a$ over $K$ in $L$, as required. This gives that the residue field $vK$ is stably embedded. A similar argument replacing the conjugation of residue field elements with value group elements gives that the value group is also stably embedded.
\end{proof}

We now consider the constant subfield. For an arbitrary differentially henselian field $(K, v, \delta)$, it is not the case that the constant subfield $(C_K, v)$ is stably embedded as a pure valued field in $K$. To see this, consider the following:

\begin{eg}[{\cite[Proposition 2.5]{CubidesDelon2016}}]
Let $(K, v) \models \ACVF_{0,0}$. Let $(K, v) \prec (L, w)$ be an elementary extension such that $K$ is dense in $L$. Let $a \in L \setminus K$, and consider the following definable set:
\[
D = \{(y, z) \in K^2 : v(a - y) > v(a - z) \}.
\]
Suppose that $D$ is definable in with parameters in $K$. Then, the family of open balls
\[
D_z = \{y \in K: (y, z) \in D\}
\]
is definable in $K$. Observe that $D_z$ is a nested family of nonempty balls with empty intersection. However, every model of $\ACVF$ is definably maximal, i.e. every definable family of nested non-empty balls has nonempty intersection (see \cite[Definition 2.1]{CubidesDelon2016}).

In particular, where $(L, w, \partial)$ is a model of $\DCVF_{(0, 0)} = \ACVF_{(0, 0} \cup \DHen$, taking $(K, v)$ to be its constant subfield equipped with the induced valuation, we have that $(K, v) \preceq (L, w)$ as pure valued fields, and $K$ is dense in $L$ with respect to the valuation topology. Applying the above, we find that $D$ is not definable with parameters in $K$.
\end{eg}

Thus, we expand the language to include sets definable in the pair of pure valued fields $(K, C_K)$. Let $\L_\vf^P$ denote the language $\L_\vf \cup \{P\}$, where $P$ is a unary predicate (later to be interpreted as a distinguished subfield). Let $K$ be a differentially henselian field, and endow $K$ with a $\L_\vf^P$-structure by interpreting $P(K) = C_K$.

Now, we will equip $C_K$ with the structure induced from the pair $(K, C_K)$. For every $\L_\vf^P(K)$-definable subset $A$ of $(C_K)^n$, we let $P_A$ be a new $n$-ary predicate. Define
\[
\L_C = \L_\vf \cup \{P_A: A \subseteq (C_K)^n \text{ is $\L_\vf^P(K)$-definable} \}.
\]
We equip $C_K$ with an $\L_C$-structure in the natural way, by interpreting $P_A$ as $A \subseteq (C_K)^n$ for each such predicate $P_A$.

\begin{lem}\label{qf_definable_sets_of_constants}
    Let $\L$ be a relational expansion of the language of rings, and let $\L^\delta = \L\cup\{\delta\}$. Let $K$ be an $\L^\delta$-structure which is a differential field. Suppose $A$ is quantifier-free definable by an $\L^\delta(K)$-formula $\phi(\bar{x})$. Then, $A \cap (C_K)^{\bar{x}} = B \cap (C_K)^{\bar{x}}$ for some quantifier-free $\L(K)$-definable set $B$.
\end{lem}
\begin{proof}
    It is sufficient to prove this statement in the case where $A$ is definable by an atomic $\L^\delta(K)$-formula and taking the relevant boolean combinations.

    Every atomic $\L^\delta(K)$-formula $\phi(\bar{x})$ is of the form $P((f_i(\bar{x}))_{i<n})$, where $P$ is an $n$-ary relation symbol of $\L$, and $f_i(\bar{x})$ is a differential polynomial in $K\{\bar{x}\}$ for each $i < n$. Let $\bar{x} = (x_0,...,x_m)$.

    Writing $f_i = (f_i)_\alg(x_0, x_0',...,x_0^{(k)},...,x_m^{(k)})$, we set 
    \[
    g_i(x_0,...,x_m) = (f_i)_\alg(x_0,0,...,0,x_1,0,...,0,x_m,0,...,0).
    \]
    That is, $g_i$ is the algebraic polynomial in $K[\bar{x}]$ obtained from $(f_i)_\alg$ by setting $x_i^{(j)} = 0$ for all $j > 0$.

    Let $\psi(\bar{x})$ be the $\L(K)$-formula $P(g_0(\bar{x}),...,g_{n-1}(\bar{x}))$, and let $B = \psi(K^{\bar{x}})$. We claim that $B \cap (C_K)^m = A \cap (C_K)^m$. We observe that for any $\bar{a} \in (C_K)^m$, we have $f_i(\bar{a}) = g_i(\bar{a})$, since every element of $\bar{a}$ is a constant. Thus, we have that
    \[
    K \models \forall \bar{x} ((\bar{x} \in (C_K)^m) \to (\phi(\bar{x}) \leftrightarrow \psi(\bar{x}))).
    \]
    Finally, we have that $A \cap (C_K)^m = B \cap (C_K)^m$, as required.
\end{proof}

\begin{thm} \label{induced_str_constants_is_pair_str_thm}
    Let $(K, v, \delta)$ be a differentially henselian field. Then, $C_K$ is stably embedded in $K$ as an $\L_C$-structure.
\end{thm}
\begin{proof}
    Let $T$ be the $\L_\vf$-theory of $K$, and let $T^\Mor$ be its Morleyisation with corresponding language $\L_\Mor$. By \cite[Theorem 2.4.2]{CubidesPoint2019}, the theory $T^\Mor \cup \mathrm{DH}$ admits quantifier elimination in the language $\L_\Mor(\delta)$. Since $K$ (with its $\L_\Mor(\delta)$-structure) is a model of $T^\Mor\cup \mathrm{DH}$, it eliminates quantifiers in the language $\L_\Mor(\delta)$.

    Let $\phi(\bar{x})$ be an $\L_\vd(K)$-formula, and let $A = \phi(K^{\bar{x}})\cap (C_K)^{\bar{x}}$. As $T^\Mor \cup \mathrm{DH}$ has quantifier elimination in $\L_\Mor(\delta)$, we have that 
    \[
    T^\Mor \cup \mathrm{DH} \models \forall\bar{x} (\phi(\bar{x}) \leftrightarrow \psi(\bar{x}))
    \]
    for some quantifier-free $\L_\Mor(\delta, K)$-formula $\psi(\bar{x})$.

    By Lemma \ref{qf_definable_sets_of_constants}, there is a quantifier-free $\L_\Mor(K)$-formula $\chi(\bar{x})$ such that $\chi(K^{\bar{x}}) \cap (C_K)^{\bar{x}} = A$. As $\chi$ is a $\L_\Mor(K)$-formula, it is equivalent modulo $T^\Mor$ to an $\L_\vf(K)$-formula $\theta(\bar{x})$, possibly with quantifiers. Thus, $A$ is $\L_\vf^P$-definable, and so by construction, $A$ is also $\L_C$-definable in $C_K$.
\end{proof}

\section{Constructing Differentially Henselian Fields} \label{section_constructions}

In this section, we will construct various examples of differentially henselian fields, namely by an inductive construction based on transcendence bases, and by iterated power series extensions.

We begin by considering a construction via iterated power series (cf. \cite[Example 5.2(ii)]{SanchezTressl2020}). Let $(K, v, \delta)$ be a valued-differential field, henselian as a pure valued field. Let $(K_0, v_0, \delta_0) = (K, v, \delta)$, and for each $n<\omega$, let 
\[ \textstyle
(K_{n+1}, v_{n+1}, \delta_{n+1}) = (K_n((t_n^\Q)), v_n \circ v_{t_n}, \hat\delta_n + \frac{\dd}{\dd t_{n}}).
\]
Define
\[
(K_\infty, v_\infty, \delta_\infty) = \bigcup_{n<\omega} (K_n, v_n, \delta_n)
\]
to be the union of the chain that we have constructed. 

% We can explicitly describe the valuation $v_\infty: K_\infty^\times \to vK \times \Q^{<\omega}$
% as follows: set $v_0: K_0 \to vK$ to be the valuation of $K$, and now suppose that we have constructed the valuation $v_n: K_n \to vK \times \Q^n$. Let $a = \sum_{\gamma \in \Q} a_\gamma t^\gamma$, and define
% \[
% v_{n+1} (a) = (v_n(a_{\gamma_0}), \gamma_0)
% \]
% where $\gamma_0 = \min(\supp(a))$. The valuation $v_\infty: K_\infty \to vK \times \Q^{<\omega}$ is then given by the union $\bigcup_{n<\omega} v_n$. 

% There is a natural valuation $v_\infty: K_\infty \to vK \times \Q^{<\omega}$, ordered reverse lexicographically, where we think of $\Q^\omega$ as the set of sequences in $\Q^\omega$ with finite support. We construct $v_\infty$ as follows:

We can explicitly describe the valuation $v_\infty$ as follows: we let $v_0$ be the valuation $v : K^\times \to vK$, and suppose we have constructed $v_n: K_n \to vK \times \Q^n$, ordered reverse lexicographically. Then, let $a = \sum_{i \in \Q} a_i t_n^i$ be an element of $K_{n+1} = K_n((t_n^\Q))$. We then define
\[
v_{n+1}(a) = (v_n(a_N), N) \in vK \times \Q^{n+1}
\]
where $N = \min\supp(a)$. Ordered reverse lexicographically, this gives precisely the valuation $v_{n+1}$ on $K_{n+1}$. Taking the union, we have that the valuation $v_\infty$ is valued in $vK \times \Q^{<\omega}$, where $\Q^{<\omega}$ is considered as the subset of elements of $\Q^\omega$ with finite support, ordered reverse lexicographically.

\begin{thm} \label{K_infty_dh}
The field $(K_\infty, v_\infty, \delta_\infty)$ is differentially henselian. 
\end{thm}
\begin{proof}
Firstly, $(K_\infty, v_\infty)$ is henselian as it is a union of henselian valued fields. We show this by induction: for $n = 0$, $v_0 = v$ is henselian by assumption. Assume $v_n$ is henselian. Then, $v_{n+1} = v_n \circ v_{t_n}$ by definition. $v_n$ and $v_{t_n}$ are both henselian, thus $v_{n+1}$ is henselian.

Now, let $f \in K_\infty\{x\}$ be a differential polynomial of order $k$, and let $\bar{a} \in K_\infty^{k+1}$ such that $f_\alg(\bar{a}) = 0$ and $s(f)_\alg(\bar{a}) \neq 0$. Let $\gamma = (q, \sigma) \in vK \times \Q^{<\omega} = v_\infty{K_\infty}$. Let $N < \omega$ such that $f \in K_N\{x\}$, $\bar{a} \in K_N$, and $N > \max(\supp(\sigma))$.

Let $T^*$ denote the standard twisted $K_N$-Taylor morphism for $K_N[[t]]$, considered as a $K_N$-subalgebra of $K_N((t^\Q))$. (We drop the subscript $N$ of $t$ for readability.) Let $A = K_N\{x\}/I(f)$. Define an algebraic $K_N$-algebra homomorphism $\phi: A \to K_N$ as follows: for $i \leq k$, set $\phi(x^{(i)}) = a_i$, and for $i > k$, define $\phi(x^{(i)})$ recursively by taking derivatives of the relation $\phi(f(x)) = 0$ and rearranging.

Consider the element $\alpha = T^*_\phi(x) \in K_{N+1}$. We claim that, for $n \leq k$, the constant term of $\alpha^{(n)}$ is precisely $a_n$. Let $\alpha = \sum_i \alpha_i t^i$, and recall that the coefficients are given by the formula
\[
\alpha_i = \frac{1}{i!} \sum_{j \leq i} (-1)^{i-j}\binom{i}{j} \partial^{i-j} (\phi(\delta^{j} x)),
\]
where $\partial$ and $\delta$ denote the derivations on $K_N$ and $A$, respectively. Thus, for $i \leq j$, we have that
\[
\alpha_i = \frac{1}{i!} \sum_{j \leq i} (-1)^{i-j}\binom{i}{j} \partial^{i-j} (a_j).
\]
Now, we compute the constant term of $(\partial + \ddt)^n (\sum_i \alpha_i t^i)$. Observe that this is precisely
\begin{align*}
\sum_{m=0}^n m! \binom{n}{m}  \partial^{n-m} \alpha_m &= \sum_{m=0}^n \binom{n}{m} \sum_{j=0}^m (-1)^{m-j}\binom{m}{j} \partial^{n-j} (a_j) \\
&= \sum_{m=0}^n \sum_{j=0}^m (-1)^{m-j} \frac{n!}{m! (n-m)!} \frac{m!}{j! (m-j)!} \partial^{n-j}(a_j)
\end{align*}
Collect the terms by the index $j$, and observe that the sum is equal to:
\[
\sum_{j=0}^n \sum_{m=j}^n (-1)^{m-j} \frac{n!}{(n-m)!} \frac{1}{j! (m-j)!} \partial^{n-j}(a_j).
\]
Setting $l = m-j$, we rewrite the sum as follows:
\begin{align*}
\sum_{j=0}^n \sum_{l=0}^{n-j} &(-1)^l \frac{n!}{(n-j-l)!j!l!} \partial^{n-j}(a_j) \\
&= \sum_{j=0}^n \sum_{l=0}^{n-j} (-1)^l \frac{n!}{j!(n-j)!} \frac{(n-j)!}{l!(n-j-l)!} \partial^{n-j}(a_j) \\
&= \sum_{j=0}^n \binom{n}{j} \sum_{l=0}^{n-j} (-1)^l \binom{n-j}{l} \partial^{n-j}(a_j).
\end{align*}
Since $\sum_{l=0}^{n-j} (-1)^l \binom{n-j}{l}$ is simply an alternating sum of binomial coefficients, it is 1 precisely when $n=j$, and 0 otherwise. We therefore conclude that the constant term of $\alpha^{(n)}$ is precisely $a_n$. 

Thus, for $i \leq k$, $v_\infty(\alpha^{(i)} - a_i) > \gamma$, as $\alpha^{(i)} - a_i$ is a series in $t_{n+1}$ with zero $t_{n+1}$-constant term, and $v(t_{n+1}) > \gamma$. Since $T^*_\phi$ is a differential $K_N$-algebra homomorphism, $\alpha$ is a solution to $f$ with $\Jet_k(\alpha) \in B_\gamma(\bar{a})$, as required.
\end{proof}

We note that, for the same reason as in the differentially large construction, it suffices to take $K((t))$ in place of $K((t^\Q))$. We therefore obtain:
\begin{cor} \label{iterated_laurent_series_dh}
Let $(K, v, \delta)$ be a valued-differential field which is henselian as a pure valued field. Then, the union of the chain 
\[ \textstyle
(K, v, \delta) \subseteq (K((t_0)), v \circ v_{t_0}, \delta + \frac{\dd}{\dd t_0}) \subseteq (K((t_0))((t_1)), v \circ v_{t_0} \circ v_{t_1}, \delta + \frac{\dd}{\dd t_0} + \frac{\dd}{\dd t_1}) \subseteq ...
\]
is differentially henselian.
\end{cor}

We recall that a field is large if and only if it is existentially closed in a henselian field. The above construction provides an analogous result for differential fields:

\begin{prop}
A differential field $(K, \delta)$ is differentially large if and only if there is a differentially henselian field $(L, w, \partial)$ with $(L, \partial) \supseteq (K, \delta)$ such that $(K, \delta)$ is existentially closed in $(L, \partial)$ as a differential field.
\end{prop}
\begin{proof}
For the backwards direction, it suffices to observe that $(L, \partial)$ is differentially large, and any differential subfield $(K, \delta)$ of $(L, \partial)$ satisfying $(K, \delta) \preceq_\exists (L, \partial)$ is also differentially large.
% We show that $(K, \delta)$ satisfies the condition in Proposition \ref{diff_large_lifting_condition}. Let $f, g \in K\{x\}$ be differential polynomials with $\ord(g) < \ord(f) = n$, and let $\bar{a} \in K$ such that $f_\alg(\bar{a}) = 0$, $g(\bar{a}) \neq 0$ and $s(f)_\alg(\bar{a}) \neq 0$. Considering $f_\alg$ and $g_\alg$ as polynomials in $L$, by continuity, there is $\gamma \in vL$ such that $B_\gamma(\bar{a})$ contains no root of $g$. Now, by differential henselianity of $(L, w, \partial)$, there is $b \in L$ such that $f(b) = 0$, and $\Jet_n(b) \in B_\gamma(\bar{a})$. In particular, we have that $g(b) \neq 0$. Finally, by existential closure of $(K, \delta)$ in $(L, \partial)$, there exists $c \in K$ such that $f(c) = 0$ and $g(c) \neq 0$.

For the forwards implication, we first note that by 4.3(iii) of \cite{SanchezTressl2020}, $(K, \delta)$ is existentially closed in $(K_n, \delta_n) \coloneqq (K((t_0))((t_1))...((t_n)), \hat\delta + \frac{\dd}{\dd t_1} + ... + \frac{\dd}{\dd t_n})$ for any $n < \omega$. Thus $(K, \delta)$ is existentially closed in the union $(L, \partial) = \bigcup_{n<\omega} (K_n, \delta_n)$. By Corollary \ref{iterated_laurent_series_dh}, $(L, w, \partial)$ is differentially henselian, as required.
\end{proof}

We will now show that any henselian valued field with sufficiently large transcendence degree admits a derivation which induces the structure of a differentially henselian field. This is an adaptation of the construction in \cite[Theorem 4.3]{SanchezTressl2023} to the henselian case.

\begin{defn}[cf. {\cite[Notation 4.1]{SanchezTressl2023}}]
Let $(K, v, \delta)$ be a valued-differential field. A \textbf{differentially henselian problem over $K$ of order $n$} is a triple $(f, \bar{c}, \gamma)$, where $f \in K\{x\}$ is a differential polynomial of order $n$, $\bar{c} \in K^{n+1}$ satisfies $f_\alg(\bar{c}) = 0$ and $s(f)_\alg(\bar{c}) \neq 0$, and $\gamma \in vK$.  We say that $(f, \bar{c}, \gamma)$ is \textbf{irreducible} if $f$ is irreducible.

Let $(L, w, \partial)$ be a valued-differential field extension of $(K, v, \delta)$. An element $a \in L$ is a \textbf{solution} to the differentially henselian problem $(f, \bar{c}, \gamma)$ if $f(a) = 0$ and $\Jet_n(a) \in B_\gamma(\bar{c})$. We say that $L$ \textbf{solves} $(f, \bar{c}, \gamma)$ if it contains a solution $a$ to $(f, \bar{c}, \gamma)$. 
\end{defn}

It is clear that a valued-differential field is differentially henselian if and only if it is henselian and solves all differentially henselian problems over itself (Theorem \ref{dh_equiv_model_DH_thm}). Further, by Proposition \ref{dh_iff_dh_for_irreducible_prop}, it suffices to solve all irreducible differentially henselian problems over itself. We begin by constructing a solution for a single differentially henselian problem.

\begin{lem} \label{singer_problems_finite}
Let $(K, v, \delta)$ be a valued-differential field, and let $(f, \bar{c}, \gamma)$ be an irreducible differentially henselian problem over $K$ of order $n$. Let $(L, w)$ be a henselian valued field extending $(K, v)$ as a pure valued field, such that $\trdeg(L/K) \geq n$. Then, there is a derivation on $L$ extending the derivation on $K$ such that there is a solution $a \in L$ to $(f, \bar{c}, \gamma)$ with $\Jet_{n-1}(a)$ algebraically independent over $K$.
\end{lem}
\begin{proof}
Since $f_\alg(\bar{c}) = 0$ and $s(f)_\alg(\bar{c}) \neq 0$, and as $s(f)_\alg(\bar{c}) \in K^\times$, by the Implicit Function Theorem for Henselian Fields (Theorem \ref{henselian_ift}), there is $\mu \in vK$ such that there is a unique continuous function $g: U = B_\mu(c_0,...,c_{n-1}) \to L$ with $f_\alg(\bar{y}, g(y)) = 0$ for any $\bar{y} \in U$. By shrinking the ball if necessary, and by continuity of $g$, we may also assume that for any $\bar{y} \in U$, $(\bar{y},g(y)) \in B_\gamma(\bar{c})$. 

We claim that $U$ contains a point of transcendence degree $n$ over $K$. Since $L/K$ is an extension of transcendence degree of at least $n$ over $K$, there is a tuple $(b_0,...,b_{n-1}) \in L^n$ algebraically independent over $K$. By replacing $b_i$ with $b_i^{-1}$ if necessary, we may also assume that $v(b_i) \geq 0$ for each $i$. Take $d \in K^\times$ with $v(d) > \mu$. Then, we claim that the tuple $(a_0,...,a_{n-1}) \coloneqq (db_0 + c_0,...,db_{n-1} + c_{n-1})$ suffices. 

If the $a_i$ are not algebraically independent over $K$, then there is a nonzero polynomial $h(x_0,..,x_{n-1}) \in K[x_0,...,x_{n-1}]$ such that $h(a_0,...,a_{n-1}) = 0$. Then, defining
\[
\tilde{h}(x_0,...,x_{n-1}) = h(dx_0 + c_0,...,dx_{n-1} + c_{n-1}),
\]
and since $d, c_0,...,c_{n-1} \in K$ with $d \neq 0$, we have that $\tilde{h} \in K[x_0,...,x_{n-1}]$ is nonzero with $\tilde{h}(b_0,...,b_{n-1}) = 0$, which contradicts the algebraic independence of the $b_i$ over $K$. Further, $v(a_i - c_i) = v(db_i) = v(d) + v(b_i) \geq \mu$ as we have assumed that $v(b_i) \geq 0$ and $v(d) \geq \mu$. Thus $(a_0,...,a_{n-1})$ is a point of transcendence degree $n$ in $U$.

Setting $a_n = g(a_0,...,a_{n-1})$ and defining $a_i' = a_{i+1}$ for $i< n$, we obtain a derivation on $K(a_0,...,a_n)$. This derivation extends uniquely to its relative algebraic closure in $L$. We now extend this arbitrarily to a derivation $\partial$ on $L$.

Now observe that, by construction, $a = a_0$ is a solution to the differentially henselian problem $(f, \bar{c}, \gamma)$ in $(L, w, \partial)$.
\end{proof}

\begin{note}
We may construct the derivation $\partial$ on $L$ such that $(L, \partial)$ is a differentially algebraic extension of $(K, \delta)$: we observe that every element of the algebraic closure of $K[a_0,...,a_n]$ in $L$ is differentially algebraic over $K$. Extending $(a_0,...,a_{n-1})$ to a transcendence basis $B$ of $L/K$, and by setting $\partial(b) = 0$ for every $b \in B\setminus\{a_0,...,a_{n-1}\}$, we obtain a derivation $\partial$ on $L$ such that every element of $L$ is differentially algebraic over $(K, \delta)$, as required.
\end{note}

We now show with an inductive construction that given a henselian valued field extension $(K, v) \subseteq (L, w)$ of sufficient transcendence degree, it is possible to construct a derivation on $L$ such that every irreducible differentially henselian problem over $K$ has a solution in $L$.

\begin{lem} \label{singer_problems_step_one}
Let $(K, v, \delta)$ be a valued-differential field, and let $(L, w)$ be a henselian valued field extension of $(K, v)$ with $\trdeg(L/K) \geq |K|$. Then, there is a derivation $\partial$ on $L$ extending $\delta$ such that $(L, w, \partial)$ solves all irreducible differentially henselian problems over $K$.
\end{lem}
\begin{proof}
Let $|K| = \kappa$, and enumerate all irreducible differentially henselian problems over $K$ by $(S_\alpha)_{\alpha < \kappa}$. This is possible, as the set of irreducible differentially henselian problems is a subset of $K\{x\} \times K^{<\omega} \times vK$ which has cardinality $\kappa$. Let $B \subseteq L$ be a transcendence basis for $L$ over $K$, and let $(B_\alpha)_{\alpha < \kappa}$ be a partition of $B$ such that $|B_\alpha| \geq \aleph_0$ for each $\alpha < \kappa$.

For $\alpha < \kappa$, let $L_\alpha$ be the relative algebraic closure of $K\left(\bigcup_{\beta < \alpha}B_\beta \right)$. In particular, $L_\alpha$ is henselian for all $\alpha < \kappa$. Set $\partial_0 = \delta$, fix some $\alpha < \kappa$, and suppose we have constructed $\partial_\beta$ for all $\beta < \alpha$ such that $(L_\beta, \delta_\beta)$ solves the differentially henselian problem $S_\mu$ for every $\mu < \beta$.

For $\alpha = \xi+1$, since $S_\xi$ is a differentially henselian problem over $L_\xi$ of finite order, by Lemma \ref{singer_problems_finite}, there is a derivation $\partial_{\alpha}$ extending $\partial_{\xi}$ such that $(L_\alpha, \partial_\alpha)$ solves $S_\xi$. When $\alpha$ is a limit, let $\partial_\alpha = \bigcup_{\beta<\alpha} \partial_\beta$.

Take $\partial = \bigcup_{\alpha < \kappa} \partial_\alpha$, and observe that $\partial$ is a derivation on $L$ such that $(L, w, \partial)$ all differentially henselian problems over $K$.
\end{proof}

\begin{note}
By applying the previous remark in every step, we may construct $\partial$ such that $(L, \partial)$ is a differentially algebraic extension of $(K, \delta)$. 
\end{note}

\begin{thm}\label{constructing_diff_hen_fields_trdeg}
Let $(K, v, \delta)$ be a valued-differential field, and $(L, w)$ be an extension of $(K, v)$ which is henselian and has $\trdeg(L/K) \geq |K|$. Then, there is a derivation $\partial$ on $L$ extending $\delta$ such that $(L, w, \partial)$ is differentially henselian.
\end{thm}
\begin{proof}
We construct a derivation $\partial$ on $L$ such that $(L, w, \partial)$ solves all irreducible differentially henselian problems over itself. Let $\kappa = \trdeg(L/K) \geq |K|$, and let $B$ be a transcendence basis for $L$ over $K$. Let $(B_n)_{n<\omega}$ be a partition of $B$ such that $|B_n| = \kappa$ for each $n < \omega$. Let $L_n$ be the relative algebraic closure of $K\left( \bigcup_{m<n} B_m\right)$ in $L$. In particular, $L_n$ is henselian for every $n < \omega$.

By construction, $\trdeg(L_{n+1}/L_n) = \kappa$ for each $n<\omega$. We construct the derivation $\partial$ inductively. Let $\partial_0$ be the unique extension of $\delta$ to the relative algebraic closure of $K$ in $L$. Suppose we have constructed $\partial_n$ on $L_n$ such that $(L_n, \partial_n)$ solves all differentially henselian problems over $L_m$ for each $m < n$. Applying Lemma \ref{singer_problems_step_one}, we find a derivation $\partial_{n+1}$ on $L_{n+1}$ such that $(L_{n+1}, \partial_{n+1})$ solves all irreducible differentially henselian problems over $L_n$. 

Take $\partial = \bigcup_{n<\omega} \partial_n$. Since every (irreducible) differentially henselian problem over $L$ is an irreducible differentially henselian problem over $L_n$ for some $n < \omega$, and since $(L, w, \partial)$ solves all irreducible differentially henselian problems over $L_n$ for every $n < \omega$, we have that $(L, w, \partial)$ solves all irreducible differentially henselian problems over itself, as required.
\end{proof}

Again by the previous remark, we may construct each derivation $\partial_n$ such that $(L_n, \partial_n) \subseteq (L_{n+1}, \partial_{n+1})$ is differentially algebraic for every $n$. Thus, we obtain the following:

\begin{cor}\label{inf_trdeg_diff_hensel_algebraic}
Let $(K, v, \delta)$ be a valued-differential field, and let $(L, w)$ be a henselian extension of $(K, v)$ with $\trdeg(L/K) \geq |K|$. Then, there is a derivation $\partial$ on $L$ extending $\delta$ such that $(L, w, \partial)$ is differentially henselian, and the extension of differential fields $(L, \partial)/(K, \delta)$ is differentially algebraic.
\end{cor}

We also obtain that every henselian valued field of infinite transcendence degree admits a derivation such that the resulting valued-differential field is differentially henselian.

\begin{cor}\label{inf_trdeg_diff_hensel}
Every henselian valued field $(K, v)$ of infinite transcendence degree admits a derivation $\delta$ such that $(K, v, \delta)$ is differentially henselian.
\end{cor}
\begin{proof}
Equip the prime subfield $\Q$ with the trivial derivation and induced valuation, and apply Theorem \ref{constructing_diff_hen_fields_trdeg}.
\end{proof}

% \begin{rmk}
%     We should remark here that if a valued-differential field $(K, v, \delta)$ is differentially henselian, then it necessarily has infinite transcendence degree: 
% \end{rmk}

We recall that given any differential field $(K, \delta)$ which is large as a pure field, there is a differential field extension $(L, \partial)$ such that $(L, \partial)$ is differentially large and $K \preceq L$ as pure fields (see \cite[Corollary 4.8(ii)]{SanchezTressl2020} and \cite[Theorem 6.2(II)]{Tressl2005}). From the above, we can extract an analogous result for valued fields, and recover \cite[Theorem 2.3.4]{CubidesPoint2019} for the differentially henselian case.

\begin{cor} \label{elt_extn_admits_dh_der}
    Let $(K, v, \delta)$ be a valued-differential field which is henselian as a pure valued field. There is a valued-differential field extension $(L, w, \partial)$ of $(K, v, \delta)$ such that $(L, w, \partial)$ is differentially henselian, and $(K, v) \preceq (L, w)$ as pure valued fields.
\end{cor}
\begin{proof}
    Let $(L, w)$ be an elementary extension of $(K, v)$ of sufficiently large transcendence degree, and apply Theorem \ref{constructing_diff_hen_fields_trdeg} to obtain a derivation $\partial$ extending $\delta$ such that $(L, w, \partial)$ is differentially henselian.
\end{proof}

\section{Sections of Differentially Henselian Fields}
\label{section_sections}

It is known that every henselian valued field $(K, v)$ of equicharacteristic 0 admits a section $Kv \to \O_v$ of the residue map, i.e. an embedding $\phi: Kv \to \O_v$ such that $\res \circ \phi = \id_{K_v}$. More precisely:

\begin{prop} \label{henselian_sections_algebraic}
Let $A$ be a local henselian ring, with maximal ideal $\mathfrak{m}$ and residue field $L$. If $\ch(L) = 0$, then $\res: A \to L$. In fact, for any maximal subfield $K \subseteq A$, the restriction $\res|_K: K \to L$ is surjective.
\end{prop}

We will show that an analogous result holds for differentially henselian fields of equicharacteristic 0, under certain saturation conditions.

\begin{prop}\label{henselian_sections_differential}
Let $(K, v, \partial)$ be a differentially henselian field of equicharacteristic 0. Then:
\begin{enumerate}[label=(\roman*)]
\item For any maximal differential subfield $L \subseteq \O_v$, the restriction $\res|_K: L \to Kv$ is surjective.
\item Suppose $K$ is $\aleph_1$-saturated. Then, for any arbitrary derivation $\delta$ on $Kv$, there exists a section $Kv \to \O_K$ of $\res$ respecting the derivation $\delta$.
\end{enumerate}
\end{prop}
\begin{proof}
For (i), let $L \subseteq \O_K$ be a maximal differential subfield such that the restriction of the residue map to $L$ is not surjective onto $Kv$. Let $\hat{L}$ be a (not necessarily differential) maximal subfield of $\O_v$ containing $L$. Let $a \in \hat{L}\setminus L$. By differential henselianity, there exists a constant $b \in \O_K$ such that $v(b-a) > 0$, i.e. $v(b) = 0$ and $\res(b) = \res(a)$. Then, $L\la b \ra = L(b)$ is a differential subfield of $\O_v$ strictly containing $L$, contradicting maximality.

For (ii), it suffices by Proposition \ref{henselian_sections_algebraic} to find a maximal subfield $L$ of $\O_v$ such that the restriction $\res|_L: L \to Kv$ is a differential isomorphism. We inductively construct such a maximal subfield. Enumerate the elements of $Kv$ as $(a_\alpha)_{\alpha < \kappa}$ for some appropriate cardinal $\kappa$.

Let $L_0 = \Q \subseteq \O_v$, and suppose we have constructed the differential subfield $L_\alpha \subseteq \O_v$ for some ordinal $\alpha < \kappa$. For convenience of notation, we will identify $L_\alpha$ with its image in $Kv$ under the residue map. We consider three cases for $a_\alpha$.

Firstly, if $a_\alpha \in L_\alpha$, then let $L_{\alpha+1} = L_\alpha$. 

Now, suppose that $a_\alpha$ is differentially algebraic over $L_\alpha$ in $(Kv, \delta)$, and let $f$ be the differential minimal polynomial of $a_\alpha$ over $L_\alpha$, and let $n = \ord(f)$. Let $F$ be any maximal (not necessarily differential) subfield of $\O_v$ extending $L_\alpha$. By Proposition \ref{henselian_sections_algebraic}, the residue map restricted to $F$ is surjective onto $Kv$. Let $\bar{a} \in F^{n+1}$ denote the preimage of $\Jet_n(a_\alpha)$. By construction, $f_\alg(\bar{a}) = 0$, and $s(f)_\alg(\bar{a}) \neq 0$. 

As $K$ is differentially henselian, there exists $b \in K$ such that $\Jet_n(b) \in B_0(\bar{a})$. In particular, we have that $\Jet_n(b) \in \O_v^{n+1}$ and for $i \leq n$, $\res(\partial^i(b)) = \res(a_i) = \delta^i(a_\alpha)$. We claim that the differential minimal polynomial of $b$ over $L_\alpha$ is $f$. 

Suppose there is a differential polynomial $g$ over $L_\alpha$ of order $m < n$ such that $g(b) = 0$. Then, since $\res(\partial^i b) = \delta^i a_\alpha$ for every $i \leq n$, and passing to the residue field, we see that $a_\alpha$ is also a root of $g$, contradicting minimality of $f$.

We set $L_{\alpha+1} = L_\alpha\la b\ra$. We claim that $L_{\alpha+1}$ is a differential subfield of $\O_v$. To see this, suppose that $c$ is some nonzero polynomial combination of $b,...,\partial^n(b)$ over $L_\alpha$. That is, for some polynomial $h$ over $L_\alpha$, $c = h(b,...,\partial^n(b))$. By construction, $c = h(a_0 + m_0,...,a_n+m_n)$ for some $m_0,...,m_n$ in the maximal ideal $\mathfrak{m}_v$. By taking a Taylor expansion of $h$ at $\bar{a}$, we find that $c = h(\bar{a}) + m$, where $m \in \mathfrak{m}_v$. We note that $h(\bar{a}) \in F^\times$, hence $v(h(\bar{a})) = 0$ and $c \in \O_v^\times$, as required.

Now we consider the case where $a_\alpha$ is differentially transcendental over $L_\alpha$. By $\aleph_1$-saturation, there is an element $b \in \O_v$ such that $\res(\Jet(b)) = \Jet(a_\alpha)$. By a similar argument to the algebraic case, $b$ cannot satisfy any differential polynomial over $L_\alpha$. Let $L_{\alpha+1}=L_\alpha\la b \ra$. As before, $L_{\alpha+1}$ is a differential subfield of $\O_v$.

For a limit ordinal $\alpha$, we set $L_\alpha = \bigcup_{\beta<\alpha} L_\beta$. Taking $L = L_\kappa$, we find that this is a maximal subfield of $\O_v$ with the desired property.
\end{proof}

\begin{rmk}
    The condition that $(K, v, \partial)$ is $\aleph_1$-saturated is only needed in the case where $(Kv, \delta)$ contains differentially transcendental elements over $\Q$. 
    In the case that $(Kv, \delta)$ is differentially algebraic over $\Q$, the above result holds without $\aleph_1$-saturation.
\end{rmk}

\section{Differentially Large Henselian Fields} 

In this section, we discuss the connections between differential largeness and differential henselianity.

We apply the characterisation of differentially henselian fields in the previous section, and use a theorem of Widawski to prove that a partial converse of Lemma \ref{diff_hens_implies_diff_large_lem} holds in the case where the field is not algebraically closed. We reproduce here the relevant theorem:

\begin{thm}[\cite{Widawski2024}] \label{etale_open_differentially_large}
Let $(K, \delta)$ be a differential field which is large as a pure field. Then, $(K, \delta)$ is differentially large if and only if, for every irreducible closed variety $V \subseteq \A^{n+1}$ of dimension $n$, and not of the form $W \times \A$ for a subvariety $W$ of $\A^n$, the set $(K^{n+1})_d \cap \Reg(V) \cap V(K)$ is dense in $\Reg(V) \cap V(K)$ with respect to the \'etale-open topology on $V(K)$, where $(K^{n+1})_d \coloneqq \{\Jet_n(a): a \in K\}$, and $\Reg(V)$ denotes the regular points of $V$.
\end{thm}

We also require the following fact about the \'etale-open topology on a henselian field:

\begin{thm}[{\cite[Theorem 6.5]{JTWY2020}}]\label{etale_open_henselian}
Let $K$ be a henselian field (i.e. a field which admits a non-trivial henselian valuation) which is not algebraically closed. Then, the \'etale-open topology coincides with the unique henselian valuation topology.
\end{thm}

Combining these two results, we obtain the following:

\begin{prop} \label{diff_large_not_ac_implies_diff_hensel_prop}
    Let $(K, v, \delta)$ be a valued-differential field such that $(K, v)$ is henselian and not algebraically closed, and $(K, \delta)$ is differentially large. Then, $(K, v, \delta)$ is differentially henselian.
\end{prop}
\begin{proof}
    Let $f$ be an irreducible differential polynomial of order $n$, and suppose we have $\bar{a} \in K^{n+1}$ such that $f_\alg(\bar{a}) = 0$ and $s(f)_\alg(\bar{a}) \neq 0$. Let $V \subseteq \A^{n+1}$ denote the variety defined by the equation $f_\alg = 0$. Then, we have that $\bar{a}$ is a regular $K$-rational point of $V$. Let $\gamma \in vK$, and consider the open set $U = B_\gamma(\bar{a}) \cap V$. By Theorem \ref{etale_open_henselian}, $U$ is also open with respect to the \'etale-open topology. Now, by Theorem \ref{etale_open_differentially_large}, as $U$ contains a regular $K$-rational point, it contains a differential regular $K$-rational point $\Jet_n(b)$, as required.
\end{proof}

\begin{cor} \label{theory_of_dh_given_by_vf_and_df_cor}
Let $(K, v, \delta)$ be a valued-differential field, not algebraically closed, such that $(K, \delta)$ is differentially large, and $(K, v)$ is henselian. Then,
\[
\Th(K, v) \cup \Th(K, \delta) \models \Th(K, v, \delta).
\]
\end{cor}
\begin{proof}
As $(K, v, \delta)$ is differentially large, not separably closed, and henselian, by Proposition \ref{diff_large_not_ac_implies_diff_hensel_prop}, $(K, v, \delta)$ is a model of $(\DL)$. Thus, we may apply \cite[Corollary 2.4.7]{CubidesPoint2019}, to obtain that the theory $\Th(K, v, \delta)$ is determined by $\Th(K, v)$. 
\end{proof}

The non-algebraic closedness condition is necessary. We exhibit a class of counterexamples in the algebraically closed case:

\begin{prop} \label{dcf_acvf_not_dh}
    Every differentially closed field $(K, \delta) \models \DCF_0$ admits a nontrivial henselian valuation $v$ such that $(K, v, \delta)$ is not differentially henselian.
\end{prop}
\begin{proof}
    Let $(K, \delta)$ be a differentially closed field, and denote by $C_K$ the field of constants of $K$. Let $a \in K \setminus C_K$. As $C_K$ is algebraically closed, $a$ is transcendental over $C_K$, and $C_K[a]$ is isomorphic to the ring of polynomials $C_K[t]$. Let $\mathfrak{p}$ be the maximal ideal $aC_K[a]$ of $C_K[a]$.
    
    By Chevalley's Extension Theorem, there is a valuation ring $\O$ of $K$ containing $C_K[a]$ such that $\m \cap C_K[a] = \mathfrak{p}$, where $\m$ is the maximal ideal of $\O$. In particular, $\O$ is a nontrivial valuation ring of $K$ containing $C_K$. Let $v$ denote the valuation on $K$ corresponding to the valuation ring $\O$. As $K$ is algebraically closed, $v$ is automatically a henselian valuation.

    Now, we observe that $v(C_K^\times) = \{0\}$, as every nonzero constant is a unit of the valuation ring. Thus, there are no nonzero constants in the (open) ball $B_0(0) = \mathfrak{m}$, and $C_K$ is not dense in $K$ with respect to the valuation $v$. As the constants of any differentially henselian field are dense with respect to the valuation topology, we conclude that $(K, v, \delta)$ is not differentially henselian.
\end{proof}

\begin{cor}
    The theory $\DCF_0 \cup \ACVF_{(0, 0)}$ is not complete.
\end{cor}
\begin{proof}
    By Corollary \ref{inf_trdeg_diff_hensel}, there exist models of $\DCF_0 \cup \ACVF_{(0,0)}$ which are differentially henselian, and by Proposition \ref{dcf_acvf_not_dh}, there are models which are not differentially henselian.
\end{proof}

\section{Equivalent Characterisations of Differential Henselianity} \label{section_equiv_chars}

In this section, we will show a number of equivalent characterisations of differentially henselian fields, in the style of Theorem 4.3 of \cite{SanchezTressl2020}. 
We begin by making a modification to the notion of an (abstract) Taylor morphism (in the sense of \cite{Ng2023}) for application to differentially henselian fields, namely that a `valued' Taylor morphism should only move points infinitesimally.

\begin{defn}
    Let $(K, v, \delta) \subseteq (L, w, \partial)$ be valued-differential fields. Let $T$ be a $(K, \delta)$-Taylor morphism for $(L, \partial)$. We say that $T$ is \textbf{valued}, if, for any differential ring $A$, ring homomorphism $\phi: A \to K$ and $a \in A$, we have that 
    \[
    w(\phi(a) - T_\phi(a)) > vK.
    \]
    Alternatively, we will say that $T^*$ is a \textbf{valued $(K, v, \delta)$-Taylor morphism for $(L, w, \partial)$}.
\end{defn}

\begin{prop} \label{standard_ttm_is_valued}
    Let $(K, v, \delta)$ be a valued-differential field. Then, the standard twisted Taylor morphism $T^*$ for $(K((t)), v\circ v_t, \hat\delta + \ddt)$ is valued.
\end{prop}
\begin{proof}
    This follows by direct computation. Let $(A, \partial)$ be a differential ring, let $\phi: A \to K$ be a ring homomorphism, and let $a \in A$. Recall that $T^*_\phi(a)$ is given by the series $\sum_i b_i t^i$, where the coefficients $b_i$ are determined by the following formula:
    \[
    b_i = \frac{1}{i!} \sum_{j \leq i} (-1)^{i-j} \binom{i}{j} \delta^{i-j}(\phi(\partial^j(a))).
    \]
    Observe that $b_0 = \phi(a)$, and thus $T^*_\phi(a) - \phi(a) = \sum_{i \geq 1} b_i t^i$. Since $v_t(T^*_\phi(a) - \phi(a)) > 0$, we have that $(v \circ v_t)(T^*_\phi(a) - \phi(a)) > vK$, as required.
    \end{proof}

\begin{lem} \label{ec_in_vttm_implies_dh}
    Let $(K, v, \delta) \subseteq (L, w, \partial)$ be valued-differential fields, where $(K, v)$ is henselian as a pure valued field. Suppose that $(K, v, \delta)$ is existentially closed in $(L, w, \partial)$ as valued-differential fields, and that $(L, w, \partial)$ admits a valued $(K, v, \delta)$-Taylor morphism. Then $(K, v, \delta)$ is differentially henselian.
\end{lem}
\begin{proof}
    Let $f(x) \in K\{x\}$ be a differential polynomial of order $n$, and suppose that $\bar{a}$ is an algebraic root of $f_\alg$ with $s(f)_\alg(\bar{a}) \neq 0$. Let $\gamma \in vK$.

    Let $A$ be the differential ring $K\{x\}/I(f)$. Define the ring homomorphism $\phi: A \to K$ by evaluation at $\bar{a}$ (as $A$ is generated as a $K$-algebra by $x,..,x^{(n)}$).

    Consider $T^*_\phi: A \to L$. As $T_\phi$ is differential, $T_\phi(x)$ is a differential root of $f$. Since $T$ is valued, we have that $v(T_\phi(x^{(i)}) - \phi(x^{(i)})) > vK$, in particular, greater than $\gamma$. Hence,
    \[
    (L, w, \partial) \models \exists x (f(x) = 0 \wedge \Jet_n(x) \in B_\gamma(\bar{a})).
    \]
    By existential closure of $(K, v, \delta)$ in $(L, w, \partial)$, we obtain $b \in K$ with $f(b) = 0$ and $\Jet_n(b) \in B_\gamma(\bar{a})$ also.
\end{proof}

For a henselian valued field $(K, v)$, we know that as $K$ is large, $K$ is existentially closed in $K((t))$ as a field. Let $v \circ v_t$ denote the composition of the $t$-adic valuation $v_t$ on $K((t))$ with the valuation $v$ on the residue field. Then, we have the following:

\begin{fact} \label{henselian_vf_ec_in_laurent_fact}
    Let $(K, v)$ be a henselian valued field. Then $(K, v)$ is existentially closed in $(K((t)), v \circ v_t)$.
\end{fact}
\begin{proof}
    This is a consequence of the proof of Theorem 5.14 of \cite{Kuhlmann2016}.
\end{proof}

\begin{thm} \label{dh_iff_ec_laurent}
    Let $(K, v, \delta)$ be a valued-differential field, henselian as a pure valued field. Then, $(K, v, \delta)$ is differentially henselian if and only if $(K, v, \delta)$ is existentially closed in $(K((t)), v \circ v_t, \hat\delta + \ddt)$.
\end{thm}
\begin{proof}
    For the forward direction, suppose $(K, v, \delta)$ is differentially henselian. By Fact \ref{henselian_vf_ec_in_laurent_fact}, $(K, v)$ is existentially closed in $(K((t)), v \circ v_t)$, we have that $(K, v, \delta)$ is existentially closed in $(K((t)), v \circ v_t, \hat\delta + \ddt)$.

    For the backwards direction, suppose that $(K, v, \delta)$ is existentially closed in $(K((t)), v \circ v_t, \hat\delta + \ddt)$ as a valued-differential field. By assumption, $(K, v)$ is nontrivially henselian. By Proposition \ref{standard_ttm_is_valued}, $(K((t)), v\circ v_t, \hat\delta + \ddt)$ admits a valued $(K, v, \delta)$-Taylor morphism. By Lemma \ref{ec_in_vttm_implies_dh}, we have that $(K, v, \delta)$ is differentially henselian, as required.
\end{proof}

Iterating this argument yields a version for iterated power series extensions:
\begin{thm} \label{dh_iff_ec_laurent_iterated}
    Let $(K, v, \delta)$ be a valued-differential field, henselian as a pure valued field. Then, $(K, v, \delta)$ is differentially henselian if and only if $(K, v, \delta)$ is existentially closed in $(L, w, \partial) = (K((t_0))...((t_{n-1})), v \circ v_{t_0} \circ ... \circ v_{t_{n}}, \hat\delta + \frac{\dd}{\dd t_0} + ... + \frac{\dd}{\dd t_{n}})$ for any $n < \omega$.
\end{thm}
\begin{proof}
    Suppose that $(K, v, \delta)$ is differentially henselian. By applying Fact \ref{henselian_vf_ec_in_laurent_fact} repeatedly, we have that $(K, v)$ is existentially closed in $(L, w)$ by transitivity of existential closure. Thus $(K, v, \delta)$ is existentially closed in $(L, w, \partial)$ by differential henselianity.

    Conversely, $(L, w, \partial)$ admits a valued $(K, v)$-Taylor morphism, as we can take the inclusion of $(K((t_0)), v \circ v_{t_0}, \hat\delta + \frac{\dd}{\dd t_0})$ in $L$. We conclude by Lemma \ref{ec_in_vttm_implies_dh}, $(K, v, \delta)$ is differentially henselian.
\end{proof}

We now generalise the characterisation of differential largeness in \cite[Theorem 4.3(iv)]{SanchezTressl2020}, which says that a differential field is differentially large if and only if it is large as a pure field, and every differentially finitely generated $K$-algebra with a $K$-rational point also has a differential $K$-rational point.

\begin{thm}\label{dh_iff_4.3iv(h)}
    Let $(K, v, \delta)$ be a valued-differential field, henselian as a pure valued field. Then, $(K, v, \delta)$ is differentially henselian if and only if for every differentially finitely generated $K$-algebra $A$ with a finite differential generating set $(a_i)_{i<n}$, and a $K$-rational point $\phi: A \to K$, there is, for  any $\gamma \in vK$, a differential $K$-rational point $\psi: A \to K$ such that $v(\phi(a_i) - \psi(a_i)) > \gamma$ for each $i < n$.
\end{thm}
\begin{proof}
    Suppose $(K, v, \delta)$ is differentially henselian. Let $A$ be a differentially finitely generated $K$-algebra, $\bar{a} = (a_i)_{i<n}$ a finite differential generating set for $A$. Suppose $A$ has a $K$-rational point $\phi: A \to K$, and let $\gamma \in vK$. 

    Since $A$ is differentially finitely generated, $A$ is isomorphic to a differential $K$-algebra of the form $K\{x_0,...,x_{n-1}\}/I$, where $I$ is the kernel of the map $\pi$ which evaluates $x_i$ at $a_i$. Write $\bar{x} = (x_0,...,x_{n-1})$.
    
    Applying the standard twisted Taylor morphism to $\phi$, we obtain a differential $K$-algebra homomorphism $T^*_\phi: A \to K((t))$. By Proposition \ref{standard_ttm_is_valued}, for any $f \in A$, we have that $(v \circ v_t)(\phi(f) - T^*_\phi(f)) > vK$. So, in particular, $T^*_\phi(\bar{a}) \in B_\gamma(\phi(\bar{a}))$.

    Let $\mathfrak{p}$ denote the preimage of $\ker(T^*_\phi)$ under $\pi$ in $K\{\bar{x}\}$. Clearly, $\mathfrak{p}$ is a differential prime ideal of $K\{\bar{x}\}$, thus by the Ritt-Raudenbush basis theorem, it is finitely generated as a differential radical ideal of $K\{\bar{x}\}$. Let $f_0(\bar{x}),...,f_{k-1}(\bar{x})$ be such a generating set.

    By construction, we now have that $T^*_\phi(\bar{a})$ is a solution to the system of differential polynomials $f_0(\bar{x}) = ... = f_{k-1}(\bar{x}) = 0$. As $(K, v, \delta)$ is differentially henselian, it is existentially closed in $(K((t)), v \circ v_t, \hat\delta + \ddt)$ by Theorem \ref{dh_iff_ec_laurent}. 

    Applying existential closure, there is some tuple $\bar{b} \in K$ with $f_i(\bar{b}) = 0$ for each $i < k$, and also $\bar{b} \in B_\gamma(\phi(\bar{a}))$. Equivalently, the differential prime ideal $\mathfrak{p}$ vanishes on $\bar{b}$, and so in particular, vanishes on $I$. 
    
    Thus, the map $K\{\bar{x}\} \to K$ defined by evaluation at $\bar{b}$ factors through $K\{\bar{x}\}/I = A$, thus $A$ has a differential point $\psi: A \to K$ with $\psi(\bar{a}) = \bar{b}$. Further, $\psi(\bar{a}) = \bar{b} \in B_\gamma(\phi(\bar{a}))$, as required.
\end{proof}

From the proof above, it is easy to see that we may reduce to the case where the differentially finitely generated $K$-algebra $A$ is a domain. Another variation on this condition is that we may require the image of the differential generators $\bar{a}$ to be `close up to order $n$' for some finite $n$. We state these more precisely as follows:

\begin{prop}
    Let $(K, v, \delta)$ be a valued-differential field, henselian as a pure valued field. The following are equivalent:
    \begin{enumerate}[label=(\roman*)]
        \item $(K, v, \delta)$ is differentially henselian.
        \item For any differentially finitely generated $K$-algebra $A$ with a finite differential generating set $\bar{a}$ and $K$-rational point $\phi: A \to K$, for any $\gamma \in vK$, there is a differential $K$-rational point $\psi: A \to K$ such that $\psi(\bar{a}) \in B_\gamma(\phi(\bar{a}))$.
        \item For any differentially finitely generated $K$-algebra $A$, which is a domain, with a finite differential generating set $\bar{a}$ and $K$-rational point $\phi: A \to K$, for any $\gamma \in vK$, there is a differential $K$-rational point $\psi: A \to K$ such that $\psi(\bar{a}) \in B_\gamma(\phi(\bar{a}))$.
        \item For any differentially finitely generated $K$-algebra $A$ with a finite differential generating set $\bar{a}$ and $K$-rational point $\phi: A \to K$, for any $\gamma \in vK$ and $n < \omega$, there is a differential $K$-rational point $\psi: A \to K$ such that $\psi(\Jet_n(\bar{a})) \in B_\gamma(\phi(\Jet_n(\bar{a})))$.
    \end{enumerate}
\end{prop}
\begin{proof}
    The equivalence of (i) and (ii) is given by Theorem \ref{dh_iff_4.3iv(h)}. The implications (ii) $\implies$ (iii) and (iv) $\implies$ (ii) are trivial. It remains to show (iii) $\implies$ (ii) and (ii) $\implies$ (iv). For (iii) $\implies$ (ii), as in the proof of Theorem \ref{dh_iff_4.3iv(h)}, we may replace $A$ with $K\{\bar{x}\}/\pi^{-1}(\ker(T^*_\phi))$ and proceed. For (ii) $\implies$ (iv), we simply observe that $\Jet_n(\bar{a})$ is also a finite differential generating set for $A$, and apply (ii).
\end{proof}

An interesting corollary of this result is that isolated $K$-rational points of differentially finitely generated $K$-algebras must be differential:

\begin{defn} \label{isolated_k_rat_pt}
    Let $(K, v, \delta)$ be a valued-differential field, and let $A$ be a differentially finitely generated $K$-algebra. We say that a $K$-rational point $\phi: A \to K$ is \textbf{isolated} (with respect to the valuation topology) if there is a finite differential generating set $\bar{a} = (a_i)_{i<n}$ and $\gamma \in vK$ such that for any $K$-rational point $\psi: A \to K$ with $\psi(\bar{a}) \in B_\gamma(\phi(\bar{a}))$, we have that $\phi = \psi$.
\end{defn}

\begin{cor}
    Let $(K, v, \delta)$ be differentially henselian, and let $A$ be a differentially finitely generated $K$-algebra. Let $\phi: A \to K$ be an isolated $K$-rational point. Then $\phi$ is differential.
\end{cor}
\begin{proof}
    As $\phi$ is isolated, there is a finite differential generating set $\bar{a} = (a_i)_{i<n}$ of $A$ and $\gamma \in vK$ such that there is no $K$-rational point $\psi: A \to K$ with $\psi \neq \psi$ and  $\psi(\bar{a}) \in B_\gamma(\phi(\bar{a})$.

    Theorem \ref{dh_iff_4.3iv(h)}, there is a differential $K$-rational point $\chi: K \to A$ with $\chi(\bar{a}) \in B_\gamma(\phi(\bar{a}))$. By the above, we necessarily have that $\phi = \chi$, thus $\phi$ is already differential.
\end{proof}

\section{Differential Weil Descent on Valued Fields} \label{section_dh_weil}

In the paper \cite{SanchezTressl2018}, L\'eon S\'anchez and Tressl use a differential version of the Weil descent to show that algebraic extensions of differentially large fields are themselves differentially large (see \cite[Theorem 5.10]{SanchezTressl2020}). In this section, we will adapt some of their machinery to prove a corresponding result for differentially henselian fields.

We first recall the relevant result for differentially large fields:

\begin{thm}[{\cite[Theorem 6.1]{SanchezTressl2018}}]\label{alg_extn_of_dl_is_dl}
    Let $(K, \delta)$ be a differentially large field, and let $(L, \partial)$ be a differential field extension where $L/K$ is algebraic. Then, $(L, \partial)$ is differentially large.
\end{thm}
In particular, the algebraic closure of any differentially large field is a model of $\DCF_0$.

% We are able to obtain a partial result in the case where the field is not real closed using previous results:

% \begin{prop} \label{alg_extn_of_not_rcf_dh_is_dh}
%     Let $(K, v, \delta)$ be a differentially henselian field, where $K$ is not real closed. Let $(L, w, \partial)$ be an algebraic valued-differential field extension of $(K, v, \delta)$. Then, $(L, w, \partial)$ is differentially henselian.
% \end{prop}
% \begin{proof}
%     It suffices to show that every subextension of finite degree is itself differentially henselian. Thus we may assume that $L$ itself is an extension of finite degree of $K$.

%     Since $(K, v)$ is henselian, the extension $w$ of $v$ to $L$ is uniquely determined. In particular, $(L, w)$ is also henselian. Further, as $K$ is not real closed, $L$ is not algebraically closed. Further, $(L, \partial)$ is differentially large by Theorem \ref{alg_extn_of_dl_is_dl}. We then conclude by Theorem \ref{diff_large_hensel_implies_diff_hensel} that $(L, w, \partial)$ is differentially henselian.
% \end{proof}

From this and Proposition \ref{diff_large_not_ac_implies_diff_hensel_prop}, we can see that an algebraic extension of any differentially henselian field which is not real closed is again differentially henselian. We will establish the unrestricted result later as Theorem \ref{alg_extn_of_dh_is_dh}. To do this, we will require the use of the differential Weil descent. We begin by constructing the classical Weil descent, closely following the setup of \cite{SanchezTressl2018}. For full details and algebraic technicalities, we direct the reader to the aforementioned paper.

Let $K$ be a ring, and $L$ be a $K$-algebra. Let
\[
F: K\Alg \to L\Alg
\]
be the \textbf{extension of scalars} functor, given by $F(A) = A \otimes_K L$ on objects, and for a morphism $\phi: A \to B$, $F(\phi) = \phi \otimes_K \id_L$. We assume that tensor products are taken over $K$, unless otherwise stated, and suppress the relevant subscripts. The \textbf{Weil descent} functor $W: L\Alg \to K\Alg$ we construct shall be the left adjoint of $F$.

We assume that $L$ is a free and finitely generated $K$-module of dimension $l$ over $K$. Fix a basis $b_1,...,b_l$ of $L$ as a $K$-module. For each $i$, define $\lambda_i: L \to K$ by
\[
\lambda_i \left( \sum_j a_j b_j \right) = a_i,
\]
i.e. $\lambda_i(x)$ is the $i$th coordinate of $x$ with respect to the basis $b_1,...,b_l$. For a $K$-algebra $A$, define $\lambda_i^A = \id_A \otimes \lambda_i : A \otimes L \to A \otimes K = A$.

\begin{defn}[{\cite[Definition 2.3]{SanchezTressl2018}}] 
    Let $T$ be a set of indeterminates. Define a $K$-algebra
    \[
    W(L[T]) = K[T]^{\otimes l}.
    \]
    For $i = 1,...,l$ and $t \in T$, write
    \[
    t(i) = 1 \otimes ... \otimes 1 \otimes t \otimes 1 \otimes ... \otimes 1 \in K[T]^{\otimes l}
    \]
    where the $t$ occurs in the $i$th position. Define the $L$-algebra homomorphism $W_{L[T]}: L[T] \to F(W(L[T])) = K[T]^{\otimes l} \otimes L$ by setting for each $t \in T$
    \[
    W_{L[T]}(t) = \sum_i (t(i) \otimes b_i).
    \]
    Define $F_{K[T]} : W(F(K[T])) = K[T]^{\otimes l} \to K[T]$ by setting, for each $t \in T$ and $i = 1,...,l$:
    \[
    F_{K[T]}(t(i)) = \lambda_i(1) t.
    \]
\end{defn}

We may choose $W(L[T])$ to be the Weil descent of $L[T]$, and $W_{L[T]}$ to be the unit of the adjunction at $L[T]$. We then obtain, for any $K$-algebra $A$
\begin{align*}
    \tau = \tau(L[T], A) : \Hom_{A\Alg}(K[T]^{\otimes l}, A) &\to \Hom_{L\Alg}(L[T], A \otimes L) \\
    \phi &\mapsto F(\phi) \circ W_{L[T]}
\end{align*}
a bijective map. We can explicitly compute, for $\phi: K[T]^{\otimes l} \to A$ and $t \in T$
\[
\tau(\phi)(t) = \sum_i \phi(t(i))\otimes b_i \in A \otimes L.
\]
Conversely, if $\psi: L[T] \to A \otimes L$ is an $L$-algebra homomorphism, we obtain the corresponding $K$-algebra homomorphism $\phi: K[T]^{\otimes l} \to A$ by setting, for each $t \in T$, $i = 1,...,l$:
\[
\phi(t(i)) = \lambda_i^A(\psi(t)).
\]
Now, we may construct the Weil descent of an arbitrary $L$-algebra $B$. Let $\pi_A: L[T] \to B$ be a surjective $L$-algebra homomorphism for some set $T$ of indeterminates. Let $I_B$ be the ideal generated in $W(L[T]) = K[T]^{\otimes l}$ by all elements of the form
\[
\lambda_i^{W(L[T])}(W_{L[T]}(f))
\]
where $i = 1,...,l$, and $f \in \ker \pi_B$. Define
\[
W(B) = W(L[T])/I_B,
\]
and set
\[
W(\pi_B) : W(L[T]) \to W(B)
\]
to be the residue map.

Then, for any $K$-algebra $A$, the bijection $\tau(L[T], A)$ from above induces a bijection
\[
\tau(B, A) : \Hom_{K\Alg}(W(B), A) \to \Hom_{L\Alg}(B, F(A))
\]
such that the following square commutes:
\[
\begin{tikzcd}
\Hom_{K\Alg}(W(B), A) \ar[d, "-\circ W(\pi_B)" ] \ar[r, "\tau(B{,}A)"] & \Hom_{L\Alg}(B, F(A)) \ar[d, "-\circ \pi_B"] \\
\Hom_{K\Alg}(W(L[T]), A) \ar[r, "\tau({L[T]}{,}A)"] & \Hom_{L\Alg}(L[T], F(A)).
\end{tikzcd}
\]
From the commutativity of the above diagram, for $\phi \in \Hom_{K\Alg}(W(B), A)$ we obtain the following:
\[
\tau(B, A)(\phi)\circ \pi_B = \tau(L[T], A)(\phi\circ W(\pi_B)) = ((\phi \circ W(\pi_B))\otimes \id_L) \circ W_{L[T]}.
\]
Explicitly computing $W_B = \tau(B, W(B))(\id_{W(B)})$, for any $t \in T$, we obtain:
\[
W_B(\pi_B(t)) = \sum_i W(\pi_B)(t(i))\otimes b_i = \sum_i (t(i) + I_B) \otimes b_i.
\]

We will now show that, where $L/K$ is a finite algebraic extension of valued fields, $B$ an $L$-algebra, and $\phi, \psi: B \to L$ are $L$-rational points, then $\phi$ and $\psi$ are `close' if the corresponding $K$-rational points of the Weil descent $W(B)$ are also `close'.

\begin{nota}
    Let $B$ be an $L$-algebra, and let $a \in B$. Let $\pi_B: L[T] \to B$ be a surjective $L$-algebra homomorphism, and assume there is $t \in T$ with $\pi_B(t) = a$. Write
    \[
    a(i) = W(\pi_B)(t(i)) \in W(B).
    \]
    Equivalently, $a(i) = \lambda^{W(B)}_j(W_B(a))$.
\end{nota}

In particular, we should note that the definition of $a(i)$ has no dependence on the choice of indeterminates $T$ or the surjective homomorphism $\pi_B$.

\begin{thm}\label{weil_descent_continuous}
    Let $(K, v)$ be a valued field, let $(L, w)$ be a finite algebraic extension of $(K, v)$, and fix $b_1,...,b_l$ a basis of $L$ as a $K$-vector space. Let $B$ be an $L$-algebra, let $a \in B$ and fix $\gamma \in wL$. Then, for any pair of $K$-algebra homomorphisms $\tilde\phi, \tilde\psi : W(B) \to K$ satisfying the property that 
    \[
    v(\tilde\phi(a(i))-\tilde\psi(a(i))) > \gamma - \epsilon,
    \]
    for each $i=1,...,l$, where $\epsilon = \min_i w(b_i)$, we have that
    \[
    w(\phi(a) - \psi(a)) > \gamma,
    \]
    where $\phi, \psi: B \to L$ denote the images of $\tilde\phi$ and $\tilde\psi$ under $\tau(B, K)$, respectively.
\end{thm}
\begin{proof}
    We compute $\phi(a)$ in terms of the $a(i)$. Let $\pi_B$ be a surjective $L$-algebra homomorphism $L[T] \to B$, where there is $t \in T$ with $\pi_B(t) = a$. Then:
    \begin{align}
        \phi(a) &= (((\tilde\phi \circ W(\pi_B))\otimes \id_L) \circ W_{L[T]})(t) \\
        &= ((\tilde\phi \circ W(\pi_B))\otimes \id_L)\left( \sum_i t(i) \otimes b_i \right) \\
        &= \sum_i (\tilde\phi \circ W(\pi_B))(t(i)) \otimes b_i \\
        &= \sum_i \tilde\phi(a(i)) \otimes b_i.
    \end{align}
    Replacing $\tilde\phi$ with $\tilde\psi$, we obtain a corresponding result for $\psi$. Now, we see that
    \[
    \phi(a) - \psi(a) = \sum_i (\tilde\phi(a(i)) - \tilde\psi(a(i))) \otimes b_i.
    \]
    Taking $w$ and applying the ultrametric inequality, we obtain
    \[
    w(\phi(a) - \psi(a)) \geq \min_i (w(\tilde\phi(a(i)) - \tilde\psi(a(i))) + w(b_i)).
    \]
    Suppose that 
    \[
    v(\tilde\phi(a(i))-\tilde\psi(a(i))) > \gamma - \epsilon
    \]
    holds for each $i$. Then,
    \begin{align}
        w(\phi(a) - \psi(a)) &\geq \min_i w(\tilde\phi(a(i)) - \tilde\psi(a(i))) + \min_i w(b_i) \\
         &> \gamma - \epsilon + \epsilon = \gamma.
    \end{align}
    This shows the desired inequality.
\end{proof}

We will now show that the converse also holds in a restricted setting, where we assume that the extension $(K, v) \in (L, w)$ admits a separated basis. For full details on separated extensions, we refer the reader to \cite{Delon1988}.

\begin{defn}
    Let $(K, v) \subseteq (L, w)$ be an extension of valued fields. We say that a finite set of elements $(b_i)_{i<n}$ is \textbf{separated} if, for any $(a_i)_{i<n} \in L^n$, we have that
    \[
    w\left( \sum_i a_i b_i \right) = \min_i w(a_i b_i).
    \]
    We say that the extension $(K, v) \subseteq (L, w)$ is \textbf{separated} if every finite dimensional $K$-vector subspace of $L$ admits a separated basis.
\end{defn}

\begin{rmk}
    Clearly, any separated set of elements is linearly independent.
\end{rmk}

\begin{egs}
    \begin{enumerate}
        \item \cite[Corollary 7]{Delon1988} Every algebraic extension $(K, v) \subseteq (L, w)$ of henselian valued fields of equicharacteristic 0 is separated. 
        \item \cite[Proof of Corollary 7]{Delon1988} Any finite algebraic extension $(K, v) \subseteq (L, w)$ with 
        \[
        [L:K] = (vL: vK) \cdot [Lw : Kv]
        \]
        is separated.
    \end{enumerate}
    
\end{egs}

\begin{prop} \label{weil_descent_continuous_converse}
    Let $(K, v) \subseteq (L, w)$ be a finite algebraic extension of valued fields. Let $b_1,...,b_l$ be a separated basis of $L$ as a $K$-vector space. Let $B$ be an $L$-algebra, and let $a \in B$. Let $\phi, \psi: B \to L$ be $L$-algebra homomorphisms, and denote their images under $\tau(B, K)^{-1}$ by $\tilde\phi$ and $\tilde\psi$, respectively. Then, for each $i$,
    \[
    v(\tilde\phi(a(i))-\tilde\psi(a(i))) \geq w(\phi(t) - \psi(t)) - w(b_i).
    \]
\end{prop}
\begin{proof}
    We claim that
    \[
    \tilde\phi(a(i)) = \lambda_i(\phi(a)).
    \]
    To see this, we define $\tilde\phi$ by setting $\tilde\phi(a(i)) = \lambda_i(\phi(a))$, and we apply $\tau(B, K)$:
    \[
    \tau(B, K)(\tilde\phi)(a) = \sum_i \tilde\phi(a(i)) \otimes b_i = \sum_i \lambda_i(\phi(a))\otimes b_i = \phi(a).
    \]
    Thus $\tilde\phi = \tau(B, K)^{-1}(\phi)$.

    Now observe the following:
    \[
        w(\phi(a)-\psi(a)) = w\left(\sum_i (\tilde\phi(a(i)) - \tilde\psi(a(i))) \otimes b_i \right).
    \]
    As the basis $(b_i)$ is separated, we have that for any $j$,
    \begin{align}
        w(\phi(a) - \psi(a)) &= \min_i (w(\tilde\phi(a(i)) - \tilde\psi(a(i))) + w(b_i)) \\
        & \leq w(\tilde\phi(a(j)) - \tilde\psi(a(j))) + w(b_j).
    \end{align}
    Rearranging, we obtain the inequality
    \[
    w(\tilde\phi(a(j)) - \tilde\psi(a(j))) \geq w(\phi(a) - \psi(a)) - w(b_j).
    \]
    as required.
\end{proof}

We now recall the main result on the differential version of the Weil descent. For a differential ring $(K, \delta)$, we denote the category of differential $(K, \delta)$-algebras by $(K, \delta)\Alg$.

\begin{thm}[{\cite[Theorem 3.4]{SanchezTressl2018}}] \label{diff_Weil} 
    Let $(K, \delta)$ be a differential ring, and $(L, \partial)$ be a differential $K$-algebra, finitely generated and free as a $K$-module. Then:
    \begin{enumerate}[label=(\roman*)]
        \item The functor $F^\mathrm{diff}: (K, \delta)\Alg \to (L, \partial)\Alg$ which sends a differential $(K, \delta)$-algebra $(A, \eta)$ to $(A \otimes L, \eta \otimes \id_L + \id_A \otimes \partial)$ has a left adjoint $W^\mathrm{diff}: (L, \partial)\Alg \to (K, \delta)\Alg$ known as the \textbf{differential Weil descent} from $(L, \partial)$ to $(K, \delta)$. The differential Weil descent sends a differential $(B, \dd)$ algebra to $(W(B), \dd^W)$, where $\dd^W$ is uniquely determined in \cite[Theorem 3.2]{SanchezTressl2018}.
        \item Let $(A, \eta) \in (K, \delta)\Alg$ and $(B, \dd) \in (L, \partial)\Alg$. Then, the bijection
        \[
        \tau(A, B): \Hom_{K\Alg}(W(B), A) \to \Hom_{L\Alg}(B, F(A))
        \]
        from the classical Weil descent restricts to a bijection
        \[
        \Hom_{(K, \delta)\Alg}(W^\mathrm{diff}(B, \dd), (A, \eta)) \to \Hom_{(L, \partial)\Alg}((B, \dd), F^\mathrm{diff}(A, \eta)).
        \]
    \end{enumerate}
\end{thm}

We apply this now to show that every finite algebraic extension of a differentially henselian field is again differentially henselian.

\begin{prop}
    Let $(K, v, \delta)$ be a differentially henselian field, and let $(L, w, \partial)$ be a finite algebraic extension of $(K, v, \delta)$. Then, $(L, w, \partial)$ is differentially henselian.
\end{prop}
\begin{proof}
    Let $f \in L\{x\}$ be a differential polynomial of order $n$, and let $\bar{a} \in L$ such that $f_\alg(\bar{a}) = 0$ and $s(f)_\alg(\bar{a}) \neq 0$. Fix $\gamma$ in $wL$. Let $(B, \dd)$ be the differential $L$-algebra $L\{x\}/I(f)$, and $\phi: B \to L$ be the $L$-algebra homomorphism given by evaluation at $\bar{a}$. 

    Observe that $B$ is generated as an $L$-algebra by $x, x',..., x^{(n-1)}$. We take $\pi_B: L[t_0,...,t_{n-1}]$ to be the morphism which sends $t_i$ to $x^{(i)}$.
    Applying the Weil descent, and the bijection $\tau(B, K)$, we obtain a $K$-algebra homomorphism $\tilde\phi: W(B) \to K$. We also have that $W(B)$ is finitely generated as a $K$-algebra by the $t_i(j)$, where $i = 0,...,{n-1}$ and $j = 1,...,l$. Let $\eta \in vK$ be an element satisfying the condition in Theorem \ref{weil_descent_continuous}.

    As $(K, v, \delta)$ is differentially henselian, it is existentially closed in $(K((t)), v\circ v_t, \hat\delta + \ddt)$ by Theorem \ref{dh_iff_ec_laurent}, and as $(W(B), \dd^W)$ is differentially finitely generated by the $t_i(j)$ via the surjective homomorphism $W(\pi_B): W(L[T]) \to K$, we may apply Theorem \ref{dh_iff_4.3iv(h)} to obtain a differential $K$-algebra homomorphism $\tilde\psi: (W(B), \dd^W) \to (K, \delta)$ such that for each $i, j$, $v(\tilde\phi(W(\pi_B)(t_i(j))) - \tilde\psi(W(\pi_B)(t_i(j)))) > \eta$.

    Write $\psi = \tau(B, K)(\tilde\psi)$, and by Theorem \ref{diff_Weil}, we have that $\psi: (B, \dd) \to (L, \partial)$ is a differential $L$-algebra homomorphism. Now, by Theorem \ref{weil_descent_continuous}, we have that $w(\phi(\pi_B(t_i)) - \psi(\pi_B(t_i))) > \gamma$ for each $i$.

    Taking $b = \psi(\pi_B(t))$, we observe that $f(b) = 0$, and $b^{(i)} = \psi(\pi_B(t_i))$. Further, by the above, $w(a_i - b_i) > \gamma$ for each $i$. Thus, $(L, w, \partial)$ is differentially henselian, as required.
\end{proof}

This extends simply to arbitrary algebraic extensions:

\begin{thm} \label{alg_extn_of_dh_is_dh}
    Let $(K, v, \delta)$ be a differentially henselian field. Then every algebraic extension $(L, w, \partial)$ of $(K, v, \delta)$ is differentially henselian.
\end{thm}
\begin{proof}
    By the previous proposition, every finite subextension of $(L, w, \partial)/(K, v, \delta)$ is differentially henselian. 
    Let $f \in L\{x\}$ be a differential polynomial of order $n$, let $\bar{a} \in L^{n+1}$ be such that $f_\alg(\bar{a}) = 0$ and $s(f)_\alg(\bar{a}) \neq 0$. Let $\gamma \in wL$. Let $F$ be the extension of $K$ generated by the coefficients of $f$ and $\bar{a}$. Then, since $F/K$ is finite, $(F, u, \dd)$ (taking the appropriate restrictions of $w$ and $\partial$) is differentially henselian. Let $\mu \in uF$ with $\mu > \gamma$. By differential henselianity, there is $b \in F$ with $f(b) = 0$ and $\Jet_n(b) \in B_\mu(\bar{a})$. In particular, $b \in L$ with $\Jet_n(b) \in B_\gamma(\bar{a})$. Thus $(L, w, \partial)$ is also differentially henselian, as required.
\end{proof}

Finally, we observe that the algebraic closure of a differentially henselian field must be a model of $\DCVF$:

\begin{cor}
    Let $(K, v, \delta)$ be a differentially henselian field of characteristic $(0, p)$, where $p$ is a prime or 0. Then, its algebraic closure $(\overline{K}, \bar{v}, \bar{\delta})$ is a model of $\DCVF_{(0, p)}$.
\end{cor}
\begin{proof}
    The algebraic closure $(\overline{K}, \bar{v})$ of $(K, v)$ is a model of $\ACVF_{(0, p)}$, and $(\overline{K}, \bar{v}, \bar{\delta})$ is differentially henselian by Theorem \ref{alg_extn_of_dh_is_dh}. 
\end{proof}

\bibliography{references1}
\bibliographystyle{halpha}

\end{document}